\documentclass[11pt,onecolumn]{article}
\usepackage{color,graphicx,xcolor,cite}
\usepackage{amsfonts,amssymb,amsmath,amsthm}
\usepackage{empheq}
\usepackage[ruled]{algorithm2e}
\usepackage{url}
\usepackage{multirow,array,booktabs}
\usepackage{threeparttable}
\usepackage{cases} %马峰又添加了这一条

\usepackage{latexsym}
\usepackage{color}
\definecolor{myblue}{rgb}{0,0,0.5}
\definecolor{mygreen}{rgb}{0,0.5,0}
\definecolor{myred}{rgb}{0.5,0,0}

\usepackage{amsmath,amsthm,amssymb}
\usepackage{graphicx}
\usepackage{hyperref}
\newcommand{\nn}{\nonumber}

\def \[{\begin{equation}}
\def \]{\end{equation}}

\newtheorem{theorem}{Theorem}[section]

\newtheorem{lemma}{Lemma}[section]

\newtheorem{remark}{Remark}[section]

\def \Q{{\mbox{${\cal Q}$}}}
\def \pd{{\mbox{\footnotesize$_{\!{P\!D}}$}}}
\def \dup{{\mbox{\footnotesize$_{\!{D\!P}}$}}}

\def \M{{\mbox{${\cal M}$}}}
\def \H{{\mbox{${\cal H}$}}}
\def \G{{\mbox{${\cal G}$}}}

\begin{document}
%\begin{CJK}{GBK}{song}

\begin{center}

{\Large \bf Extensions of ADMM for Separable Convex Optimization Problems with Linear Equality or Inequality Constraints}\\

\bigskip
\medskip

  {\bf Bingsheng He}\footnote{\parbox[t]{16cm}{
 Department of Mathematics,  Nanjing University, Nanjing, China.
  This author was supported by the NSFC Grant 11871029. Email: hebma@nju.edu.cn}}
 \quad
 {\bf Shengjie Xu}\footnote{\parbox[t]{16cm}{
 Department of Mathematics, Harbin Institute of Technology, Harbin, China,  and Department of Mathematics,   Southern  University of Science and Technology, Shenzhen, China. Email: xsjnsu@163.com
  }}
  \quad
 {\bf Xiaoming Yuan}\footnote{\parbox[t]{16cm}{
 Department of Mathematics, The University of Hong Kong, Hong Kong. This author was supported  by the General Research Fund from Hong Kong Research Grant Council: 12302318. Email:  xmyuan@hku.hk
  }}

\medskip

% \today

\end{center}

{\small

\parbox{0.95\hsize}{

\hrule

\medskip

{\bf Abstract.} The alternating direction method of multipliers (ADMM) proposed by Glowinski and Marrocco is a benchmark algorithm for two-block separable convex optimization problems with linear equality constraints. It has been modified, specified, and generalized from various perspectives to tackle more concrete or complicated application problems. Despite its versatility and phenomenal popularity, it remains unknown whether or not the ADMM can be extended to separable convex optimization problems with linear inequality constraints. In this paper, we lay down the foundation of how to extend the ADMM to two-block and multiple-block (more than two blocks) separable convex optimization problems with linear inequality constraints. From a high-level and methodological perspective, we propose a unified framework of algorithmic design and a roadmap for convergence analysis in the context of variational inequalities, based on which it is possible to design a series of concrete ADMM-based algorithms with provable convergence in the prediction-correction structure. The proposed algorithmic framework and roadmap for convergence analysis are eligible to various convex optimization problems with different degrees of separability, in which both linear equality and linear inequality constraints can be included. The analysis is comprehensive yet can be presented by elementary mathematics, and hence generically understandable.

\medskip

 \medskip

  \hrule

  }}

\bigskip

\section{Introduction}  \label{Sec-1-Introduction}

The alternating direction method of multipliers (ADMM) was proposed originally in \cite{GM} by Glowinski and Marrocco for solving nonlinear elliptic problems, and it has become a benchmark algorithm for solving various convex optimization problems with linear equality constraints and separable objective functions without coupled variables. Methodologically, it can be regarded as a splitting version of the classic augmented Lagrangian method (ALM) proposed in \cite{Hes,Powell}. It has found applications in an extremely broad range of areas, particularly in fields related to data science such as machine learning, computer vision, and distributed/centralized optimization. When a concrete application is considered, the original ADMM may need to be modified or specified appropriately from various perspectives to better capture the underlying structures and properties of the specific model. Some such examples include its linearized/proximal versions as proposed in \cite{He02MP,HMY-COAP}. It has also inspired many more generalized versions for solving more complicated problems, among which are a series of ADMM variants for solving multiple-block separable convex optimization problems whose objective functions consist of more than two blocks of components without coupled variables. We refer to, e.g. \cite{Boyd,Glow12,Eck12,He20}, for some survey papers about the ADMM, among a large volume of literatures. Despite the versatility and phenomenal popularity of ADMM, it remains unknown whether or not it can be extended to separable convex optimization problems with linear inequality constraints, even for two-block separable convex optimization problems.

Let us start with the canonical two-block separable convex optimization problem with linear equality constraints
\[  \label{Problem-EqC}  \min\big\{\theta_1(x) + \theta_2(y)    \;|\;  Ax+By=b,  x\in {\cal X},  y\in {\cal Y}  \big\}, \]
where $\theta_i: {\Re}^{n_i}\to {\Re} \;(i=1,2)$ are closed
proper convex functions and they are not necessarily smooth; ${\cal
X}\subseteq \Re^{n_1}$ and $ {\cal Y}\subseteq \Re^{n_2}$ are closed convex sets;
$A\in \Re^{m\times n_1}$ and $B\in \Re^{m\times n_2}$ are given matrices; $b\in \Re^m$ is a given vector. Let $\lambda \in \Re^m$ be the Lagrange multiplier and consider the Lagrangian function of the problem  (\ref{Problem-EqC})
\[   \label{LagrangianE}
      L_E(x,y,\lambda)=  \theta_1(x)  + \theta_2(y)  - \lambda^T (Ax + By  -b ), \quad   (x,y,\lambda)\in {\cal X}\times{\cal Y}\times \Re^m.
  \]
Then, the iterative scheme of the ADMM for solving (\ref{Problem-EqC}) reads as
   \[ \label{Eq-ADMM}
\hbox{(ADMM)}\quad \left\{
 \begin{array}{l}
      x^{k+1} \in \arg\min \bigl\{L_E(x, y^k,\lambda^k) + \frac{\beta}{2} \|Ax+By^k-b\|^2  \;|\;  x\in{\cal X}  \bigr\},\\[0.15cm]
    y^{k+1} \in \arg\min \bigl\{L_E(x^{k+1}, y,\lambda^k) + \frac{\beta}{2} \|Ax^{k+1}+By-b\|^2   \;|\;   y\in{\cal Y}\bigr\},\\[0.15cm]
      {\lambda}^{k+1} = \arg\max \bigl\{L_E(x^{k+1}, y^{k+1},\lambda) - \frac{1}{2\beta} \|\lambda-\lambda^k\|^2   \;|\;    \lambda\in{\Re^m}  \bigr\},
\end{array} \right.
\]
where $\beta>0$ is the penalty parameter. That is, the ADMM (\ref{Eq-ADMM}) generates the new output $(x^{k+1},y^{k+1},\lambda^{k+1})$ with the input $(y^k, \lambda^k)$. Note that the update of $\lambda^{k+1}$ in \eqref{Eq-ADMM} can be explicitly expressed as
$$
     {\lambda}^{k+1} = \lambda^k - \beta (A x^{k+1} + By^{k+1} -b).
    $$
We refer to, e.g., \cite{Gabay,Glow84,Glow89,HY-SINUM,HY-NM,YZZ-JMLR}, for some convergence study of the ADMM (\ref{Eq-ADMM}). The ADMM (\ref{Eq-ADMM}) updates the variables $x$ and $y$ by treating the functions $\theta_1$ and $\theta_2$ separately in its iterations, and the subproblems in (\ref{Eq-ADMM}) are usually much easier than the original problem (\ref{Problem-EqC}). For many application problems, the subproblems in (\ref{Eq-ADMM}) could be easy enough to have closed-form solutions or be solved up to high precisions. This feature mainly accounts for the versatility and efficiency of the ADMM (\ref{Eq-ADMM}) in various areas. Certainly, how difficult the resulting subproblems in (\ref{Eq-ADMM}) are still depends on the corresponding functions, coefficient matrices, and constraint sets; and many variants of the ADMM (\ref{Eq-ADMM}) have been proposed in the literatures for more meticulous studies. But we only concentrate on the foundational case (\ref{Eq-ADMM}), and for succinctness, we do not further elaborate on more detailed cases such as how to solve the resulting subproblems.

If the linear equality constraints in (\ref{Problem-EqC}) are changed to linear inequality constraints while all the other settings are remained, we obtain the following model:
\[ \label{Problem-I}
       \min \big\{\theta_1(x) + \theta_2(y)   \;|\;  A x  + By \ge b ,  x\in  {\cal X}, y\in {\cal Y}\big\}.
                \]
The two-block separable convex optimization model (\ref{Problem-I}) with linear inequality constraints captures particular applications such as the support vector machine with a linear kernel in \cite{cher2007,Vapn1995} and its variants in \cite{Lee2001,Man2000}. To solve (\ref{Problem-I}), it is easy to see that it can be reformulated as the following three-block separable model with linear equality constraints:
\[   \label{Problem-Eq-z}
       \min \big\{\theta_1(x) + \theta_2(y)    \;|\;   A x  + By   - z= b,    x\in  {\cal X}, y\in {\cal Y}, z\in \Re^m_+\big\},
                   \]
where $z\in \Re^m_+$ is an auxiliary variable. Then, a direct extension of the ADMM (\ref{Eq-ADMM}) can be applied to the reformulated model (\ref{Problem-Eq-z}). More specifically, let  $\lambda \in \Re^m$ be the Lagrange multiplier and the Lagrangian function of (\ref{Problem-Eq-z}) be
\[   \label{LagrangianEz}
      L_{E}(x,y,z,\lambda)=  \theta_1(x)  + \theta_2(y)  - \lambda^T (Ax + By  -z-b ), \quad   (x,y,z,\lambda)\in {\cal X}\times{\cal Y}\times \Re^m_+\times \Re^m.
  \]
Directly extending the ADMM (\ref{Eq-ADMM}) to (\ref{Problem-Eq-z}) results in the scheme
   \[ \label{Eq-ADMM-z}
\hbox{(EADMM)}\quad \left\{
 \begin{array}{l}
      x^{k+1} \in \arg\min \bigl\{L_{E}(x, y^k,z^k,\lambda^k) + \frac{\beta}{2} \|Ax+By^k-z^k-b\|^2   \;|\;    x\in{\cal X}  \bigr\},\\[0.1cm]
    y^{k+1} \in \arg\min \bigl\{L_{E}(x^{k+1}, y,z^k,\lambda^k) + \frac{\beta}{2} \|Ax^{k+1}+By-z^k-b\|^2   \;|\;    y\in{\cal Y} \bigr\},\\[0.1cm]
    z^{k+1} \in \arg\min \bigl\{L_{E}(x^{k+1}, y^{k+1},z,\lambda^k) + \frac{\beta}{2} \|Ax^{k+1}+By^{k+1}-z-b\|^2   \;|\;    z\in{\Re^m_+}\bigr\},\\[0.1cm]
      {\lambda}^{k+1} = \arg\max \bigl\{L_{E}(x^{k+1}, y^{k+1},z^{k+1},\lambda) - \frac{1}{2\beta} \|\lambda-\lambda^k\|^2   \;|\;   \lambda\in{\Re^m}  \bigr\},
\end{array} \right.
\]
where $\beta>0$ is also the penalty parameter. According to \cite{CHYY}, however, convergence of the direct extension of ADMM (\ref{Eq-ADMM-z}) is not guaranteed unless more restrictive conditions on the objective functions, coefficient matrices, as well as the penalty parameter, are additionally posed. Alternatively, the scheme (\ref{Eq-ADMM-z}) should be revised appropriately to render the convergence. For example, the output of (\ref{Eq-ADMM-z}) should be further corrected by those correction steps studied in \cite{HTY-SIOPT,HTY-MOR,HY-COAP}.

The number of blocks really matters for extensions of the ADMM (\ref{Eq-ADMM}), from both theoretical and numerical perspectives. In the literatures, there are numerous numerical studies showing that, when the ADMM (\ref{Eq-ADMM}) and its variants are applied, it is generally not preferred to artificially create more blocks of variables/functions by introducing auxiliary variables. One reason is the mentioned lack of theoretical guarantee of convergence as rigorously analyzed in \cite{CHYY}. Another more subtle reason is that if the underlying augmented Lagrangian function is decomposed by more than twice (which is usually for the sake of obtaining subproblems as easy as those in (\ref{Eq-ADMM})), then the approximation to the underlying augmented Lagrangian function might be too inaccurate and accordingly the resulting scheme may become numerically slower or even divergent. One more consequence is that tuning the penalty parameter $\beta$ usually becomes more challenging when the underlying augmented Lagrangian function is decomposed into too many blocks. This consequence is certainly based on experience and empirical study, instead of rigorous theory. Hence, for the generic two-block convex optimization model with linear inequality constraints (\ref{Problem-I}), it is interesting to discuss whether or not we can extend the original ADMM (\ref{Eq-ADMM}) in some senses that all the major features and structures of the original ADMM (\ref{Eq-ADMM}) can be kept. That is, the underlying augmented Lagrangian function should be decomposed only twice at each iteration, the resulting subproblems should be as easy as those in (\ref{Eq-ADMM}), and the convergence can be rigorously guaranteed without extra conditions. To the best of our knowledge, this question remains unknown and our main purpose is to answer this question.

Let us combine both (\ref{Problem-EqC}) and (\ref{Problem-I}) in our discussion, and consider the general two-block separable convex optimization model with linear equality or inequality constraints
\[ \label{Problem-IneC}
       \min \big\{\theta_1(x) + \theta_2(y)  \;|\;  A x  + By=b\ (\hbox{or} \ge b) ,  x\in  {\cal X}, y\in {\cal Y} \big\},
                \]
in which the settings are the same as those in (\ref{Problem-EqC}) and (\ref{Problem-I}). The solution set of the model (\ref{Problem-IneC}) is assumed to be nonempty. Our main interest is certainly the case of (\ref{Problem-IneC}) with linear inequality constraints, i.e., (\ref{Problem-I}). The reason for considering (\ref{Problem-IneC}) is that the algorithmic framework and the roadmap for convergence analysis to be presented are eligible to both (\ref{Problem-EqC}) and (\ref{Problem-I}). Another reason is that, as mentioned, we prefer to keep the features and structures of the original ADMM (\ref{Eq-ADMM}) when the linear inequality constraints are considered in (\ref{Problem-IneC}) because of both theoretical and numerical purposes. Hence, treating (\ref{Problem-EqC}) and (\ref{Problem-I}) uniformly can help us discern the difference of the to-be-proposed algorithms from the original ADMM (\ref{Eq-ADMM}) more clearly. From a high-level and methodological perspective, we will propose a unified framework of algorithmic design and convergence analysis for the model (\ref{Problem-IneC}), with which a series of specific algorithms can be easily designed and their convergence can be proved uniformly by following a common roadmap without any additional conditions on the functions, coefficient matrices, or the penalty parameter. We aim at laying down the foundation of algorithmic design and convergence analysis for extensions of the ADMM (\ref{Eq-ADMM}) from the canonical two-block model (\ref{Problem-EqC}) to the more general one (\ref{Problem-IneC}), as well as to the even more complicated multiple-block one (\ref{Problem-m}).

The rest of this paper is organized as follows. In Section \ref{Sec2-VI-TYKJ}, we review the variational inequality characterization of the model (\ref{Problem-IneC}). Our analysis will be mainly conducted in the variational inequality context. Then, we extend the ADMM (\ref{Eq-ADMM}) and propose an prototypical algorithmic framework in Section \ref{sec:algframe}; its convergence is also proved in this section. In Sections \ref{Sec3-PD} and \ref{Sec4-DP}, we specify the prototypical algorithmic framework as two concrete algorithms for the model (\ref{Problem-IneC}). In Sections \ref{Sec5-MultiB}-\ref{Sec7}, we consider a multiple-block generalized model of (\ref{Problem-IneC}), i.e., (\ref{Problem-m}), and parallelize the analysis in Sections \ref{sec:algframe}-\ref{Sec4-DP}, respectively. In Section \ref{Sec:overview}, we give an overview of how the proposed algorithmic frameworks can be unified for convex optimization problems with different degrees of separability. Finally, some conclusions are drawn in Section \ref{Sec-Conclusion}.

\section{Variational inequality characterization} \label{Sec2-VI-TYKJ}

\setcounter{equation}{0}

In this section, we summarize some preliminaries for further analysis. In particular, the variational inequality (VI) characterizations of the optimization problems appearing in our discussion are crucial. As analyzed in our previous works such as \cite{HY-SINUM,HY-COAP}, the VI approach appears to be a convenient and powerful tool for us to look into the structure of the problem under discussion, as well as to conduct convergence analysis. Our analysis in this paper will also be conducted in the context of variational inequalities.

We start from the VI representation of the optimality condition of a convex optimization problem. The following lemma will be frequently used in our following analysis. Its proof is elementary and it can be found in, e.g., \cite{Beck}.

\begin{lemma} \label{CP-TF}
\begin{subequations} \label{CP-TF0}
Consider the optimization problem
$$
\min\big\{\theta(z) + f(z)  \;|\;  z\in {\cal Z}\big\},
$$
where ${\cal Z}\subset \Re^n$ is a closed convex set, $\theta(z)$ and $f(z)$ are convex functions.
   If $f$ is differentiable on an open set which contains ${\cal Z}$, and the solution set of this  problem is nonempty, then we have that
     \[  \label{CP-TF1}   z^*  \in \arg\min \big\{ \theta(z) + f(z)   \;|\;   z\in {\cal Z}\big\}
     \]
if and only if
\[\label{CP-TF2}
      z^*\in {\cal Z}, \quad   \theta(z) - \theta(z^*) + (z-z^*)^T\nabla f(z^*) \ge 0, \quad \forall\, z\in {\cal Z}.
      \]
\end{subequations}
\end{lemma}

Now, let us focus on the model (\ref{Problem-IneC}) and derive its optimality condition in terms of the VI formulation. Without ambiguity of notation, let us reuse $\lambda$ for the Lagrange multiplier and consider the Lagrangian function of the problem  (\ref{Problem-IneC})
\[   \label{LagrangianF}
      L(x,y,\lambda)=  \theta_1(x)  + \theta_2(y)  - \lambda^T (Ax + By  -b ), \quad   (x,y,\lambda)\in {\cal X}\times{\cal Y}
        \times \Lambda,
  \]
where
   $$    \Lambda  =\left\{ \begin{array}{ll}
               \Re^m,           &   \hbox{if   $Ax+ By = b$} ,   \\[0.1cm]
                \Re^m_+,    &        \hbox{if   $Ax + By \ge b$}.
                \end{array}        \right.
     $$
Furthermore, let $\Omega:={\cal X}\times{\cal Y}\times\Lambda$. We call a point $(x^*,y^*,\lambda^*)$ defined on $\Omega$ a saddle point of the Lagrangian function (\ref{LagrangianF}) if it satisfies the inequalities
 $$    L_{\lambda\in\Lambda}(x^*,y^*,\lambda) \le L(x^*,y^*,\lambda^*) \le L_{x\in
 {\cal X}, y\in {\cal Y}}(x, y,\lambda^*). $$
Obviously, a saddle point can be characterized by
 $$
  (x^*, y^*, \lambda^*)\in \Omega, \quad    \left\{ \begin{array}{rl}
      L(x,y^*,\lambda^*)   -L(x^*,y^*,\lambda^*) \ge 0,  & \forall\, x\in {\cal X},  \\[0.1cm]
        L(x^*,y,\lambda^*)   -L(x^*,y^*,\lambda^*) \ge 0,  & \forall\, y\in {\cal Y},  \\[0.1cm]
           L(x^*,y^*,\lambda^*)  - L(x^*,y^*,\lambda)  \ge 0,       &  \forall \,\lambda\in \Lambda.
        \end{array} \right.
         $$
Alternatively, according to Lemma \ref{CP-TF}, the inequalities above can be written as the following VIs:
\[   \label{VI-Chara}
   (x^*, y^*, \lambda^*)\in \Omega, \quad   \left\{ \begin{array}{rl}
   \theta_1(x) -  \theta_1(x^*) + (x-x^*)^T(- {A}^T\lambda^*) \ge 0, & \forall\, x\in {\cal X},  \\[0.1cm]
  \theta_2(y) -  \theta_2(y^*) + (y-y^*)^T(- B^T\lambda^*) \ge 0,  & \forall\, y\in {\cal Y},  \\[0.1cm]
 (\lambda-\lambda^*)^T( Ax^* + By^* -b)\ge 0,  &  \forall \,\lambda\in \Lambda.
        \end{array} \right.
       \]
More compactly, \eqref{VI-Chara} can be rewritten as
\begin{subequations}  \label{VI-ID}
\[  \label{VI-I2}
  w^*\in \Omega, \quad   \theta(u)  - \theta(u^*) + (w-w^*)^TF(w^*) \ge 0, \quad \forall \,  w\in
     \Omega,
     \]
where
\[ \label{VI-Iw}
     w = \left(\begin{array}{c}
                      x\\
                      y\\[0.1cm]
                  \lambda \end{array} \right),
  \;\;  u = \left(\begin{array}{c}
                      x\\
                      y\end{array} \right), \;\; \theta(u)=\theta_1(x) +\theta_2(y), \;\;
    F(w) =\left(\begin{array}{c}
     - {A}^T\lambda \\
       - {B}^T\lambda \\
     {A}x + By -b \end{array} \right). \]
     \end{subequations}
It is clear that the function $\theta(u)$ is convex and  the operator $F$  in \eqref{VI-Iw} is affine with a skew-symmetric matrix. Thus, we have
\[ \label{Skew-S} (w-\tilde{w})^T(F(w)-F(\tilde{w}))\equiv 0 ,  \;\; \forall\, w, \tilde{w}.   \]
 The  solution set of the VI \eqref{VI-ID} is denoted by $\Omega^*$,  which is also  the set of the saddle points of the Lagrangian function  \eqref{LagrangianF} defined on $\Omega$.

\section{Prototypical algorithmic framework}\label{sec:algframe}

\setcounter{equation}{0}

In this section, we focus on the VI reformulation (\ref{VI-ID}), and propose an algorithmic framework conceptually in the context of variational inequalities. This algorithmic framework will be the prototype for various specific algorithms and we will show how to specify the prototypical algorithmic framework as concrete algorithms for the model (\ref{Problem-IneC}). We shall also prove the convergence of the algorithmic framework, and establish a roadmap for the convergence proof. With this roadmap, proving the convergence for different algorithms specified from the prototypical algorithmic framework can simply be reduced to identifying two matrices and then verifying the positive definiteness of another matrix. This prototypical algorithm framework and its roadmap for convergence analysis can enable us to treat a series of different algorithms uniformly from a high-level perspective, and help us discern their respective difference from the original ADMM (\ref{Eq-ADMM}) clearly.

\subsection{Algorithmic framework}\label{sec:algframe-1}

First of all, for any $w=(x, y, \lambda)\in \Re^{n_1}\times \Re^{n_2}\times \Re^m$, we define $\xi\in\Re^{3m\times3m}$ by
\[ \label{XiPw}    \xi:= Pw,  \qquad \hbox{where} \qquad
  P = \left(\begin{array}{ccc}
          \sqrt{\beta} A  &   0       & 0   \\
                0 & \sqrt{\beta} B &       0     \\
                        0   &    0    & \frac{1}{\sqrt{\beta}} I_m
   \end{array} \right) . \]
Accordingly, we define
$${\Xi}= \big\{\xi  \;|\; \xi= Pw,\, w\in \Omega\big\}\qquad   \hbox{and}  \qquad   {\Xi^*}= \big\{\xi^* \;|\; \xi^*= Pw^*,\, w^*\in \Omega^*\big\}. $$
Then, we propose a prototypical algorithmic framework for the VI reformulation (\ref{VI-ID}), which is in a prediction-correction structure.

 \begin{center}\fbox{
 \begin{minipage}{15.5cm}  {\bf{A Prototypical Algorithmic Framework for VI $\boldsymbol{(\ref{VI-ID})}$}}.

\begin{subequations}\label{M-PRECOR}
   \begin{enumerate}
   \item (Prediction Step) With $\xi^k \in \Xi$, find $\tilde{w}^k \in \Omega$ such that
  \[  \label{M-PRE} \tilde{w}^k \in \Omega, \;\;\theta(u) - \theta(\tilde{u}^k)  +  (w - \tilde{w}^k)^T F(\tilde{w}^k) \ge
      (\xi  -\tilde{\xi}^k)^T{\Q}(\xi^k-\tilde{\xi}^k), \;\; \forall \,  w \in {\Omega},
        \]
  with ${\Q} \in \Re^{3m\times 3m}$, and the matrix ${\Q}^T+{\Q}$ is positive definite.

\item (Correction Step) With $\tilde w^k$ solved by (\ref{M-PRE}) and thus ${\tilde \xi}^k:= P{\tilde w}^k$, generate $\xi^{k+1}$ by
\[ \label{M-COR}   {\xi}^{k+1} = {\xi}^k -  {\M}(\xi^k - \tilde{\xi}^k),  \]
where ${\M} \in \Re^{3m\times 3m}$ is non-singular.
\end{enumerate}
\end{subequations}
\end{minipage}} \end{center}

\subsection{Roadmap for convergence analysis}\label{sec:algframe-2}

Now, we prove the convergence of the prototypical algorithmic framework \eqref{M-PRECOR}. This is the unified analysis of convergence for various algorithms that can be specified from the prototypical algorithmic framework (\ref{M-PRECOR}). The roadmap for convergence analysis will become clear based on the following analysis.

 \begin{theorem} \label{HauptA}
 For the matrices  ${\Q}$ in \eqref{M-PRE} and $\M$ in \eqref{M-COR},  if  there is a positive definite matrix  $\H\in \Re^{3m\times3m}$  such that
  \[ \label{M-HMQ}   \H\M={\Q}  \]
  and
       \[ \label{M-HMG}   \G:=  {\Q}^T  + {\Q} - \M^T\H\M\succ 0,  \]
   then we have
  \[  \label{HauptA0}   \|{\xi}^{k+1} -{\xi}^*\|_{\H}^2
       \le   \|{\xi}^k -{\xi}^*\|_{\H}^2  -   \|{\xi}^k - \tilde{\xi}^k  \|_{\G}^2,
           \quad \forall \, \xi^*\in {\Xi^*}.  \]
   \end{theorem}
\noindent{\bf Proof}. Setting  $w$ in \eqref{M-PRE} as any fixed $w^*\in \Omega^*$, and using \eqref{Skew-S}
$$
(\tilde{w}^k-w^*)^T F(\tilde{w}^k)\equiv(\tilde{w}^k-w^*)^T F(w^*),
$$ we get
$$    (\tilde{\xi}^k-\xi^*)^T{\Q}(\xi^k-\tilde{\xi}^k)
 \ge   \theta(\tilde{u}^k)-\theta(u^*) +   (\tilde{w}^k-w^*)^T F(w^*), \quad \forall \,w^*\in \Omega^*.  $$
The right-hand side of the last inequality is non-negative. Thus, we have
    \[   \label{M-XiQ}     ({\xi}^k-\xi^*)^T{\Q}(\xi^k-\tilde{\xi}^k)  \ge   (\xi^k-\tilde{\xi}^k)^T {\Q}
       (\xi^k-\tilde{\xi}^k), \quad \forall \, \xi^*\in \Xi^*.    \]
Then, by simple manipulations, we obtain
\begin{eqnarray*}
  \lefteqn{ \|{\xi}^k -{\xi}^*\|_{\H}^2  -   \|{\xi}^{k+1} -{\xi}^*\|_{\H}^2 } \nn \\
     &\stackrel{\eqref{M-COR}}{=} & \|{\xi}^k -{\xi}^*\|_{\H}^2  -   \|({\xi}^{k} -{\xi}^*) -\M(\xi^k - \tilde{\xi}^k)  \|_{\H}^2  \nn \\
 & \stackrel{\eqref{M-HMQ}}{=}&  2 ({\xi}^{k} -{\xi}^*) ^T{\Q}(\xi^k - \tilde{\xi}^k)
       -\|\M(\xi^k - \tilde{\xi}^k)  \|_{\H}^2   \nn \\
       &\stackrel{\eqref{M-XiQ}} {\ge}  &  2 ({\xi}^{k} -\tilde{\xi}^k)^T {\Q}(\xi^k - \tilde{\xi}^k)
       -\| \M(\xi^k - \tilde{\xi}^k)  \|_{\H}^2   \nn \\
  & = & (\xi^k - \tilde{\xi}^k)^T[(  {\Q}^T +  {\Q})   -
      \M^T\H\M]  (\xi^k - \tilde{\xi}^k) \nn  \\
  &  \stackrel{\eqref{M-HMG}}{=} &  \|{\xi}^k - \tilde{\xi}^k  \|_{\G}^2.
  \end{eqnarray*}
The assertion of this theorem is proved. \qquad {$\Box$}

 \begin{theorem} \label{HauptB}  Let $\{\xi^k\}$ and $\{\tilde{w}^k\}$ be the sequences generated by the algorithmic
 framework \eqref{M-PRECOR}.
  If  the conditions \eqref{M-HMQ}  and \eqref{M-HMG} are satisfied, then the sequence $\{ \xi^k\}$ converges to  some
  $\xi^{\infty} \in \Xi^*$.
    \end{theorem}
 \noindent{\bf Proof}.  First of all,  it follows  from  \eqref{HauptA0} that the sequence  $\{\xi^k\}$
   is bounded and
   \[   \label{Key-Eq} \lim_{k\to\infty}  \|{\xi}^k - \tilde{\xi}^k  \|_{\G}^2 =0.
              \]
Thus, the sequence $\{\tilde{\xi}^k\}$ is also  bounded.  Let $\xi^{\infty}$ be a  cluster point of $\{\tilde{\xi}^k\}$ and
 $\{\tilde{\xi}^{k_j}\}$ be a subsequence converging to $\xi^{\infty}$. Let  $\{\tilde{\xi}^k\}$ and $\{\tilde{\xi}^{k_j}\}$
 be the induced sequences by $\{\tilde{w}^k\}$ and $\{\tilde{w}^{k_j}\}$, respectively. It follows from  \eqref{M-PRE} that
 $$    \tilde{w}^{k_j}\in \Omega,  \;\; \theta(u) - \theta(\tilde{u}^{k_j})  +  (w - \tilde{w}^{k_j})^T F(\tilde{w}^{k_j}) \ge
      (\xi  -\tilde{\xi}^{k_j})^T{\Q}(\xi^{k_j}-\tilde{\xi}^{k_j}), \quad \forall  \, w \in {\Omega}. $$
 Since the matrix ${\Q}$ is non-singular, it follows from \eqref{Key-Eq}, the continuity  of $\theta(u)$, and $F(w)$ that
 $$  w^{\infty}\in \Omega,  \;\; \theta(u) - \theta(u^{\infty})  +  (w - w^{\infty})^T F(w^{\infty}) \ge
    0, \quad \forall  \, w \in {\Omega}.   $$
This VI indicates that $w^{\infty}$ is a solution point of \eqref{VI-ID}, and thus $\xi^\infty=Pw^\infty\in\Xi^\ast$. Moreover, it follows from  \eqref{Key-Eq} and $\lim_{j\to\infty} \tilde{\xi}^{k_j} =\xi^{\infty}$ that the subsequence $\{\xi^{k_j}\}$ also converges to $\xi^{\infty}$.
 Finally, because of \eqref{HauptA0}, we have
 $$     \|\xi^{k+1} - \xi^{\infty} \|_{\H}^2 \le \|\xi^k - \xi^{\infty}\|_{\H}^2,   $$
 and thus $\{\xi^k\}$ converges to $\xi^{\infty}$.  The proof is complete. \qquad {$\Box$}

The convergence of the prototypical algorithmic framework (\ref{M-PRECOR}) is thus proved in Theorem \ref{HauptB}. As just shown, the proof essentially requires to verify the conditions (\ref{M-HMQ}) and (\ref{M-HMG}).

\subsection{Remarks}

Based on the analysis above, concrete algorithms for the model (\ref{Problem-IneC}) can be constructed by choosing specific matrices $\Q$ in (\ref{M-PRE}) and $\M$ in (\ref{M-COR}), and then their convergence can be proved by simply verifying the conditions (\ref{M-HMQ}) and (\ref{M-HMG}). The proposed prototypical algorithmic framework (\ref{M-PRECOR}) opens a door to designing various specific algorithms that are tailored for specific applications of the model (\ref{Problem-IneC}), and the proposed roadmap for convergence analysis provides a unified and simplified way to prove the convergence of various algorithms. Here, we give the prototypical algorithmic framework (\ref{M-PRECOR}) and the roadmap for convergence analysis from a methodological perspective, rather than discussing how to choose the matrices $\Q$ and $\M$ optimally, which should vary from case to case when a specific application is under consideration. In Sections \ref{Sec3-PD} and \ref{Sec4-DP}, we will present some such algorithms and illustrate how to follow the proposed roadmap for convergence analysis. For each algorithm, we will still use the same letters to denote these matrices but with some subscripts.

\section{Primal-dual extension of the ADMM (\ref{Eq-ADMM}) for (\ref{Problem-IneC})} \label{Sec3-PD}

\setcounter{equation}{0}

First of all, let us revisit the original ADMM (\ref{Eq-ADMM}) for the model (\ref{Problem-EqC}), and introduce an auxiliary notation
 $$     {\lambda}^{k+\frac{1}{2}}  := \lambda^k -  \beta
    (Ax^k+By^k-b).  $$
Then, ignoring some constant terms in the objective functions of the corresponding subproblems, we can rewrite the ADMM \eqref{Eq-ADMM} as
\[ \label{Eq-ADMM-A}
  \left\{\begin{array}{l}
x^{k+1}    \in \arg\!\min
 \bigl\{\theta_1(x)  - x^TA^T\lambda^{k+\frac{1}{2}} +{\textstyle{\frac{\beta}{2}}}\|A(x-x^k)\|^2  \;|\;  {x\in {\cal X}} \bigr\}, \\[0.15cm]
 y^{k+1} \in \arg\!\min\bigl\{  \theta_2(y)  - y^TB^T\lambda^{k+\frac{1}{2}}+{\textstyle{\frac{\beta}{2}}} \| B(y-y^k)+A(x^{k+1}-x^k)\|^2  \;|\;  y\in {\cal Y}  \bigr\},  \\[0.15cm]
  {\lambda}^{k+1} =  \arg\!\max  \bigl\{ - \lambda^T\bigl(Ax^{k+1}  +B{y}^{k+1}-b\bigr)  -  {\textstyle{\frac{1}{2\beta}}}\| \lambda -\lambda^k\|^2   \;|\;  \lambda\in \Re^m\bigr\}.
  \end{array} \right.
  \]
This is a reformulation of the ADMM (\ref{Eq-ADMM}) with some terms that are meticulously regrouped. It will be the reference for us to discern the difference of the new algorithms from the original ADMM (\ref{Eq-ADMM}) more conveniently.

\subsection{Algorithm}

The first algorithm for (\ref{Problem-IneC}) is presented below. Since the primal variables $x$ and $y$ are updated first before the dual variable $\lambda$, it is called a primal-dual extension of the ADMM (\ref{Eq-ADMM}) for \eqref{Problem-IneC}.

 \begin{center}\fbox{
 \begin{minipage}{15.5cm}  {\bf{ A Primal-Dual Extension of the ADMM (\ref{Eq-ADMM}) for (\ref{Problem-IneC})}}.

\begin{subequations}\label{PD-2}
   \begin{enumerate}
   \item (Prediction Step) With given $(Ax^k, By^k, \lambda^k)$, find $\tilde{w}^k=(\tilde{x}^k, \tilde{y}^k,\tilde{\lambda}^k)$ via
  \[   \label{PD-A}
  \left\{ \begin{array}{l}
    \tilde{x}^k \in \hbox{argmin}\bigl\{   \theta_1(x) -  x^TA^T{\lambda}^k  +  \frac{1}{2}\beta \|A(x-x^k)\|^2  \;|\; x\in {\cal X}\bigr\},  \\[0.2cm]
    \tilde{y}^k \in  \hbox{argmin}\bigl\{   \theta_2(y) -  y^TB^T{\lambda}^k  +  \frac{1}{2}\beta \| B(y-y^k) + A(\tilde{x}^k-x^k)\|^2  \;|\;   y\in {\cal Y}\bigr\},    \\[0.2cm]
    \tilde{\lambda}^k= \arg\!\max  \bigl\{-\lambda^T\bigl(A\tilde{x}^k +B\tilde{y}^k-b\bigr)  -     \frac{1}{2\beta}\| \lambda -\lambda^k\|^2  \;|\;  \lambda\in {\Lambda}\bigr\}.
  \end{array}  \right.
\]
\item (Correction Step) Correct the predictor $\tilde{w}^k$ solved by (\ref{PD-A}), and generate the new iterate $(Ax^{k+1}, By^{k+1}, \lambda^{k+1})$ with $\nu\in (0,1)$ by
\[ \label{PD-COR}
\left(\begin{array}{c}
Ax^{k+1} \\[0.1cm]
By^{k+1} \\[0.1cm]
\lambda^{k+1}
\end{array}\right)   = \left(\begin{array}{c}
Ax^{k} \\[0.1cm]
By^{k} \\[0.1cm]
\lambda^{k}
\end{array}\right)  -\left(\begin{array}{ccc}
            {\nu}I_m  &     -{\nu} I_m &  0   \\[0.1cm]
                 0 &   {\nu} I_m & 0  \\[0.1cm]
                - \nu\beta I_m &   0  &   I_m
        \end{array}\right) \left(\begin{array}{c}
Ax^{k} -A\tilde{x}^{k}  \\[0.1cm]
By^{k} - B \tilde{y}^{k}\\[0.1cm]
\lambda^{k} -\tilde{\lambda}^k
\end{array}\right).
\]
\end{enumerate}
\end{subequations}
\end{minipage}} \end{center}

\begin{remark}\label{PD-2-Remark 1}
Comparing with the reformulated iterative scheme of the ADMM (\ref{Eq-ADMM-A}), we see that the only difference in the prediction step (\ref{PD-A}) is the constant vector $\lambda^k$ in the crossing terms (equivalently, constant vectors in the corresponding quadratic terms), while all major features and structures of the ADMM (\ref{Eq-ADMM-A}) are remained in (\ref{PD-A}). This very minor difference does not essentially change the difficulty of the resulting $x$- and $y$-subproblems. That is, the $x$- and $y$-subproblems in (\ref{PD-A}) are of the same difficulty as those in the original ADMM (\ref{Eq-ADMM-A}) (i.e., (\ref{Eq-ADMM})). For the $\lambda$-subproblem in (\ref{PD-A}), it can be specified respectively as
$$  \tilde{\lambda}^k = \lambda^k -\beta (A\tilde{x}^k + B\tilde{y}^k -b) \qquad\hbox{or}
     \qquad  \tilde{\lambda}^k = [\lambda^k -\beta (A\tilde{x}^k + B\tilde{y}^k -b)]_+,
     $$
when the model (\ref{Problem-EqC}) or (\ref{Problem-I}) is considered, either of which is easy to compute.

\end{remark}

\begin{remark}\label{PD-2-Remark 2}
The correction step (\ref{PD-COR}) requires extremely simple computation. Looking into the implementation of the ADMM (\ref{Eq-ADMM}), we know that it is the sequence $\{Ax^k, By^k,\lambda^k\}$, instead of $\{x^k, y^k,\lambda^k\}$, that are essentially required for executing the iterations. Hence, when the ADMM (\ref{Eq-ADMM}) and its variants are implemented, one advantage is that $Ax^k$ (rather than $x^k$) and $By^k$ (rather than $y^k$) can be treated together for recursions. The correction step (\ref{PD-COR}) exactly has this advantage too, and it treats $Ax^k$, $A{\tilde x}^k$, $Bx^k$ and $B{\tilde y}^k$ aggregately with very cheap computation for updating them. Indeed, only few floating-point additions are required. Overall, the algorithm (\ref{PD-2}) maintains all major structures and features of the original ADMM (\ref{Eq-ADMM}); the resulting subproblems are of the same difficulty; and the additional computation required by the correction step (\ref{PD-COR}) is ignorable.

\end{remark}

\subsection{Specification of the prototype algorithmic framework (\ref{M-PRECOR})}

Now, we show that the algorithm (\ref{PD-2}) can be obtained by specifying the prototype algorithmic framework (\ref{M-PRECOR}). That is, we identify the specific matrices $\Q$ in (\ref{M-PRE}) and $\M$ in (\ref{M-COR}) such that (\ref{M-PRE}) and (\ref{M-COR}) can be reduced to the prediction step (\ref{PD-A}) and the correction step (\ref{PD-COR}), respectively. The specified matrices corresponding to the algorithm (\ref{PD-2}) are denoted by $\Q_{\pd}$ and ${\M}_{\pd}$, respectively. Accordingly, we divide the discussion into two subsections.

\subsubsection{Analysis for the prediction step (\ref{PD-A})}  \label{Sec31-PD}

According to Lemma \ref{CP-TF}, the predictor generated by \eqref{PD-A} satisfies $\tilde{w}^k\in \Omega$ and
$$  %\label{PD-B}
\left\{ \begin{array}{l}         \theta_1(x) - \theta_1(\tilde{x}^k)  + (x - \tilde{x}^k)^T\{-A^T\lambda^k  +\beta A^T\!A(\tilde{x}^k-x^k)\} \ge 0,    \;\; \forall \, x\in {\cal X},\\[0.1cm]
       \theta_2(y) - \theta_2(\tilde{y}^k)  + (y - \tilde{y}^k)^T\{-B^T\lambda^k  +\beta B^T\!A(\tilde{x}^k-x^k)+ \beta B^T\!B(\tilde{y}^k-y^k) \} \ge 0,     \;\; \forall \,  y\in {\cal Y}, \\[0.1cm]
  \hskip 2.8cm
   (\lambda - \tilde{\lambda}^k)^T\{\frac{1}{\beta} (\tilde{\lambda}^k-\lambda^k) + (A\tilde{x}^k + B\tilde{y}^k-b) \}\ge 0, \;\;  \forall \, \lambda\in {\Lambda}.
              \end{array}\right.
    $$
Using the VI form \eqref{VI-ID}, we can rewrite it as
\[  \label{PD-C}
     \left\{  \begin{array}{l}
        \theta_1(x) - \theta_1(\tilde{x}^k)  + (x - \tilde{x}^k)^T\{\underline{-A^T\tilde{\lambda}^k}     \\  \hskip 4.7cm + \beta A^T\!A(\tilde{x}^k-x^k)  + A^T\!(\tilde{\lambda}^k-\lambda^k)\} \ge 0,    \;\; \forall \, x\in {\cal X},  \\[0.1cm]
       \theta_2(y) - \theta_2(\tilde{y}^k)  + (y - \tilde{y}^k)^T\{\underline{-B^T\tilde{\lambda}^k}
                                    +  \beta B^TA(\tilde{x}^k-x^k)     \\
 \hskip 4.7cm +\beta B^T\!B(\tilde{y}^k-y^k)  + B^T\!(\tilde{\lambda}^k-\lambda^k)\} \ge 0,    \;\; \forall \, y\in {\cal Y},   \\[0.1cm]
       \hskip 2.5cm
   (\lambda - \tilde{\lambda}^k)^T\{(\underline{A\tilde{x}^k + B\tilde{y}^k-b})  +
              (1/\beta) \; (\tilde{\lambda}^k-\lambda^k)\}\ge 0, \;\;  \forall \, \lambda\in {\Lambda}.
  \end{array} \right.
  \]
It is easy to verify that the sum of the three underlining terms in  \eqref{PD-C} is precisely $F(\tilde{w}^k)$, where $F(\cdot)$ is defined in (\ref{VI-Iw}). Hence, we have the following lemma.

\begin{lemma}  With the  given
 $(Ax^k, By^k, \lambda^k)$,  the predictor $\tilde{w}^k\in \Omega$ generated by  \eqref{PD-A} satisfies
\begin{subequations}  \label{PD-PreD}
\[ \label{PD-PreD-F }
  \tilde{w}^k\in \Omega, \;\;\;   \theta(u) - \theta(\tilde{u}^k) +  (w- \tilde{w}^k )^T \{ F(\tilde{w}^k)
      +  Q_{\pd}(\tilde{w}^k-w^k) \} \ge 0, \quad  \forall  \, w\in\Omega,\]
where
     \[ \label{PD-Q}
     Q_{\pd}=  \left(\begin{array}{ccc}
          \beta A^TA  &   0                  &        A^T    \\
          \beta B^TA  & \beta B^TB &         B^T       \\
                        0   &    0    & \frac{1}{\beta} I_m
   \end{array} \right).
    \]
    \end{subequations}
\end{lemma}
\noindent{\bf Proof}. The assertion directly follows from  \eqref{PD-C}. \qquad {$\Box$}

\smallskip

Using the notation in \eqref{XiPw}, we can rewrite the matrix  $Q_{\pd}$ in  \eqref{PD-Q} as
\[  \label{CalQ-PD}   Q_{\pd}= P^T {\Q}_{\pd} P
  \qquad \hbox{where} \qquad  {\Q}_{\pd} = \left(\begin{array}{ccc}
             I_m &   0       &   I_m  \\
              I_m   &   I_m    &  I_m   \\
            0   &    0   & I_m
   \end{array} \right).
  \]
Note that the matrices  $Q_{\pd}\in \Re^{(n_1+n_2+m)\times(n_1+n_2+m)}$ and ${\Q_{\pd}}\in  \Re^{3m\times 3m}$ are different. Recall the notation in \eqref{XiPw}. It follows from \eqref{PD-PreD} that
\[ \label{PD-PreF}
  \tilde{w}^k\in \Omega, \;\;\;   \theta(u) - \theta(\tilde{u}^k) +  (w- \tilde{w}^k )^T F(\tilde{w}^k)
        \ge (\xi -\tilde{\xi}^k)^T{\Q}_{\pd}(\xi^k - \tilde{\xi}^k),  \quad \forall \, w\in\Omega.\]
Thus, the prediction step (\ref{PD-A}) can be specified by the prototypical prediction step (\ref{M-PRE}) with $\Q:= {\Q}_{\pd}$ as defined in (\ref{CalQ-PD}).

%For any proper positive definite matrix $H$,  we have
%\[   \label{PD-XiD} \big\langle   H ({\xi}^k -{\xi}^* ),  -H^{-1} {\Q}_{\pd}(\xi^k-\tilde{\xi}^k)  \big\rangle
%      \le -  \frac{1}{2}  ({\xi}^k -\tilde{\xi}^k)^T({\Q}_{\pd}^T + {\Q}_{\pd}) (\xi^k-\tilde{\xi}^k). \]
%Because   $({\Q}_{\pd}^T + {\Q}_{\pd}) $ is positive definite, and  $H ({\xi}^k -{\xi}^* )$ is the gradient of the unknown distance function $\frac{1}{2}\|\xi - \xi^*\|_{\H}^2$  at $\xi^k$,  namely,
%$$     \nabla \bigl(\frac{1}{2} \|\xi - \xi^*\|_{\H}^2\bigr) \Big|_{\xi =\xi^k} =  H ({\xi}^k -{\xi}^* ),  $$
% the inequality \eqref{PD-XiD} means that $-H^{-1} {\Q}_{\pd}(\xi^k-\tilde{\xi}^k)$ is a descent direction
%of the unknown distance function $\frac{1}{2}\|\xi - \xi^*\|_{\H}^2$  at $\xi^k$.

\subsubsection{Analysis for the correction step (\ref{PD-COR})}   \label{Sec32-PD}

Left-multiplying the matrix $\hbox{diag}(\sqrt{\beta}I_m,\sqrt{\beta}I_m, \frac{1}{\sqrt{\beta}}I_m)$  to both sides of  the correction step \eqref{PD-COR},  we get
$$ % \label{PD-COR-a}
\left(\begin{array}{c}
   \sqrt{\beta} Ax^{k+1} \\[0.1cm]
  \sqrt{\beta} By^{k+1} \\[0.1cm]
    \frac{1}{\sqrt{\beta}}\lambda^{k+1}
\end{array}\right)   = \left(\begin{array}{c}
\sqrt{\beta}Ax^{k} \\[0.1cm]
\sqrt{\beta}By^{k} \\[0.1cm]
      \frac{1}{\sqrt{\beta}} \lambda^{k}
\end{array}\right)  -\left(\begin{array}{ccc}
         \sqrt{\beta}   {\nu}I_m  &   -  \sqrt{\beta}{\nu} I_m &  0   \\[0.1cm]
                 0 &  \sqrt{\beta} {\nu} I_m & 0  \\[0.1cm]
                - \nu\sqrt{\beta} I_m &   0  &   \frac{1}{\sqrt{\beta}}  I_m
        \end{array}\right) \left(\begin{array}{c}
Ax^{k} -A\tilde{x}^{k}  \\[0.1cm]
By^{k} - B \tilde{y}^{k}\\[0.1cm]
\lambda^{k} -\tilde{\lambda}^k
\end{array}\right).
  $$
Then, we have
\[ \label{PD-COR-B}
\left(\begin{array}{c}
   \sqrt{\beta} Ax^{k+1} \\[0.1cm]
  \sqrt{\beta} By^{k+1} \\[0.1cm]
    \frac{1}{\sqrt{\beta}}\lambda^{k+1}
\end{array}\right)   = \left(\begin{array}{c}
\sqrt{\beta}Ax^{k} \\[0.1cm]
\sqrt{\beta}By^{k} \\[0.1cm]
      \frac{1}{\sqrt{\beta}} \lambda^{k}
\end{array}\right)  -\left(\begin{array}{ccc}
            {\nu}I_m  &   -  {\nu} I_m &  0   \\[0.1cm]
                 0 &   {\nu} I_m & 0  \\[0.1cm]
                - \nu I_m &   0  &    I_m
        \end{array}\right) \left(\begin{array}{c}
  \sqrt{\beta}(Ax^{k} -A\tilde{x}^{k})  \\[0.1cm]
 \sqrt{\beta} (By^{k} - B \tilde{y}^{k})\\[0.1cm]
    \frac{1}{\sqrt{\beta}} (\lambda^{k} -\tilde{\lambda}^k)
\end{array}\right).
\]
It follows from \eqref{XiPw} that \eqref{PD-COR-B} can be written as
\begin{subequations}  \label{PD-CorB}
\[  \label{PD-CORw}    {\xi}^{k+1} = {\xi}^k -  {\M}_{\pd}(\xi^k - \tilde{\xi}^k),  \]
where
\[  \label{M-I-PD}
   {\M}_{\pd} =
           \left(\begin{array}{ccc}
            {\nu}I_m  &     -{\nu} I_m &  0   \\[0.2cm]
                 0 &   {\nu} I_m & 0  \\[0.2cm]
                - \nu I_m &   0  &   I_m
        \end{array}\right) \quad \hbox{with} \quad \nu\in (0,1). \]
\end{subequations}
Thus, the correction step (\ref{PD-COR}) can be specified by the prototypical correction step (\ref{M-COR}) with $\M:= {\M}_{\pd}$ as defined in (\ref{M-I-PD}).

\subsection{Convergence}

Then, according to the roadmap presented in Section \ref{sec:algframe-2}, proving the convergence of the algorithm (\ref{PD-2}) can be reduced to verifying the conditions (\ref{M-HMQ}) and (\ref{M-HMG}) with the specified matrices $\Q_{\pd}$ and ${\M}_{\pd}$. Let us first identify the corresponding matrix ${\H}_{\pd}$ with the specified matrices $\Q_{\pd}$ and ${\M}_{\pd}$.

\begin{lemma}  \label{PD-H}  For any $\nu\in (0,1)$,  the matrix
\[    \label{MatrixH-PD} {\H}_{\pd} = \left(\begin{array}{ccc}
        (1+\frac{1}{\nu}) I_m  &    (1+\frac{1}{\nu}) I_m   &  I_m \\[0.2cm]
             (1+\frac{1}{\nu}) I_m  &  (1+ \frac{2}{\nu})   I_m   &  I_m  \\[0.2cm]
              I_m   &   I_m   &   I_m
        \end{array}\right)  \]
 is positive definite, and it holds that
\[   \label{PD-HMQ}   {\H}_{\pd} {\M}_{\pd}= {\Q}_{\pd} ,
   \]
where  ${\M}_{\pd}$ and $ {\Q}_{\pd}$  are defined in \eqref{M-I-PD} and \eqref{CalQ-PD}, respectively.
\end{lemma}
\noindent{\bf Proof}.  It is easy to see that ${\H}_{\pd}$ is positive definite. In addition, we have
\begin{eqnarray*}
{\H}_{\pd} {\M}_{\pd}
     &=&\left(\begin{array}{ccc}
        (1+\frac{1}{\nu}) I_m  &    (1+\frac{1}{\nu}) I_m   &  I_m \\[0.2cm]
             (1+\frac{1}{\nu}) I_m  &  (1+ \frac{2}{\nu})   I_m   &  I_m  \\[0.2cm]
              I_m   &   I_m   &   I_m
        \end{array}\right)
           \left(\begin{array}{ccc}
            {\nu}I_m  &     -{\nu} I_m &  0   \\[0.2cm]
                 0 &   {\nu} I_m & 0  \\[0.2cm]
                - \nu I_m &   0  &   I_m
        \end{array}\right)   \nn   \\
        & =&    \left(\begin{array}{ccc}
             I_m &   0       &   I_m  \\
              I_m   &   I_m    &  I_m   \\
            0   &    0   & I_m
   \end{array} \right)     \;=\;  {\Q}_{\pd}.
    \end{eqnarray*}
The assertion \eqref{PD-HMQ} is proved. \qquad {$\Box$}

\medskip

\begin{lemma}  \label{PD-MG}  For the matrices ${\Q}_{\pd}$, $ {\M}_{\pd}$  and ${\H}_{\pd}$
  defined in  \eqref{CalQ-PD}, \eqref{M-I-PD}  and \eqref{MatrixH-PD},  respectively,
the matrix
 \[ \label{PD-G}  {\G}_{\pd}:= ({\Q}_{\pd}^T +  {\Q}_{\pd})   -
      {\M}_{\pd}^T{\H}_{\pd}{\M}_{\pd}     \]
is positive definite.
\end{lemma}
\noindent{\bf Proof}. First, by elementary matrix multiplications, we have
\begin{eqnarray*}
   {\M}_{\pd}^T{\H}_{\pd}{\M}_{\pd}     &  = &     {\M}_{\pd}^T{\Q}_{\pd} =
        {\Q}_{\pd}^T  {\M}_{\pd}  \nn \\
        & =&  \left(\begin{array}{ccc}
             I_m &    I_m     &   0 \\
              0    &    I_m      &   0  \\
            I_m   &    I_m   & I_m
   \end{array} \right)  \left(\begin{array}{ccc}
            {\nu}I_m  &     -{\nu} I_m &  0   \\[0.1cm]
                 0 &   {\nu} I_m & 0  \\[0.1cm]
                - \nu I_m &   0  &   I_m
        \end{array}\right)    = \left(\begin{array}{ccc}
            {\nu}I_m  &    0  &  0   \\[0.1cm]
                 0 &   {\nu} I_m & 0  \\[0.1cm]
                0 &   0  &   I_m
        \end{array}\right).
\end{eqnarray*}
Then, it follows that
\begin{eqnarray*}
   {\cal G}_{\pd}  &= & (  {\Q}_{\pd}^T +  {\Q}_{\pd})   -
      {\M}_{\pd}^T{\H}_{\pd}{\M}_{\pd}     \nn \\
        & = &  \left(\begin{array}{ccc}
            2 I_m &    I_m     &   I_m \\
              I_m    &   2 I_m      &   I_m \\
            I_m   &    I_m   &  2I_m
   \end{array} \right)     - \left(\begin{array}{ccc}
            {\nu}I_m  &    0  &  0   \\[0.1cm]
                 0 &   {\nu} I_m & 0  \\[0.1cm]
                0 &   0  &   I_m
        \end{array}\right)= \left(\begin{array}{ccc}
            (2-{\nu})I_m  &     I_m  &  I_m   \\[0.1cm]
             I_m &   (2-{\nu}) I_m &  I_m  \\[0.1cm]
                  I_m &   I_m  &   I_m
        \end{array}\right).
\end{eqnarray*}
Thus, the matrix $ {\G}_{\pd}$ is positive definite for any $\nu\in(0,1)$.    \qquad {$\Box$}

\medskip

For the matrices ${\Q}_{\pd}$ in \eqref{CalQ-PD} and ${\M}_{\pd}$  in \eqref{PD-CorB}, the assertions in  Lemmas \ref{PD-H}  and  \ref{PD-MG} hold. Consequently, we have the following theorem which essentially implies the convergence of the algorithm (\ref{PD-2}). The proof of this theorem follows directly from Theorems \ref{HauptA} and \ref{HauptB}.

\begin{theorem}
Let $\{ \xi^k\}$ be the sequence generated by the algorithm (\ref{PD-2}). Then, we have
   \[  \label{Haupt-PD}  \|{\xi}^{k+1} -{\xi}^*\|_{ {\H}_{\pd}}^2
       \le   \|{\xi}^k -{\xi}^*\|_{ {\H}_{\pd}}^2  -   \|{\xi}^k - \tilde{\xi}^k  \|_{{\G}_{\pd}}^2,
           \quad \forall \, \xi^*\in {\Xi^*},
       \]
and thus the sequence $\{ \xi^k\}$  converges to  some  $\xi^* \in \Xi^*$.
 \end{theorem}

\section{Dual-primal extension of ADMM (\ref{Eq-ADMM}) for (\ref{Problem-IneC})} \label{Sec4-DP}

\setcounter{equation}{0}

In this section, we extend the ADMM (\ref{Eq-ADMM}) in the way that the dual variable $\lambda$ is updated first and then the primal variables $x$ and $y$ are updated. We show that this new algorithm can also be obtained by specifying the prototypical algorithmic framework (\ref{M-PRECOR}). Hence, we also follow the roadmap in Section \ref{sec:algframe-2} to prove its convergence.

\subsection{Algorithm}

The second algorithm for (\ref{Problem-IneC}) is presented below; it is called a dual-primal extension of the ADMM (\ref{Eq-ADMM}).

 \begin{center}\fbox{
 \begin{minipage}{15.5cm}  {\bf{ A Dual-Primal Extension of the ADMM (\ref{Eq-ADMM}) for (\ref{Problem-IneC})}}.

\begin{subequations}\label{DP-A}
   \begin{enumerate}
   \item (Prediction Step) With given $(Ax^k, By^k, \lambda^k)$, find $\tilde{w}^k=(\tilde{x}^k, \tilde{y}^k,\tilde{\lambda}^k)$ via
  \[ \label{DP-Pre}
\left\{\begin{array}{l}
                         \tilde{\lambda}^k= \arg\!\max \bigl\{  -\lambda^T \bigl(Ax^k +By^k-b\bigr) - \frac{1}{2\beta}\| \lambda -\lambda^k\|^2   \;|\;  \lambda\in {\Lambda}\bigr\} , \\[0.2cm]
    \tilde{x}^k \in \hbox{argmin}\bigl\{   \theta_1(x) -  x^TA^T\tilde{\lambda}^k
                         +  \frac{1}{2}\beta \|A(x-x^k)\|^2
                          \;|\;  x\in {\cal X}\bigr\},  \\[0.2cm]
      \tilde{y}^k \in \hbox{argmin}\bigl\{   \theta_2(y) -  y^TB^T\tilde{\lambda}^k
                         +  \frac{1}{2}\beta \| B(y-y^k) + A(\tilde{x}^k-x^k)\|^2
                           \;|\;  y\in {\cal Y}\bigr\}.
  \end{array} \right.
  \]

\item (Correction Step) Correct the predictor $\tilde{w}^k$ generated by (\ref{DP-Pre}), and generate the new iterate $(Ax^{k+1}, By^{k+1}, \lambda^{k+1})$ with $\nu \in (0,1)$ by
\[ \label{DP-COR}
\left(\begin{array}{c}
Ax^{k+1} \\[0.1cm]
By^{k+1} \\[0.1cm]
\lambda^{k+1}
\end{array}\right)   = \left(\begin{array}{c}
Ax^{k} \\[0.1cm]
By^{k} \\[0.1cm]
\lambda^{k}
\end{array}\right)  -\left(\begin{array}{ccc}
            {\nu}I_m  &     -{\nu} I_m &  0   \\[0.1cm]
                 0 &   {\nu} I_m & 0  \\[0.1cm]
                - \beta I_m &    - \beta I_m  &   I_m
        \end{array}\right) \left(\begin{array}{c}
Ax^{k} -A\tilde{x}^{k}  \\[0.1cm]
By^{k} - B \tilde{y}^{k}\\[0.1cm]
\lambda^{k} -\tilde{\lambda}^k
\end{array}\right).
    \]
\end{enumerate}
\end{subequations}
\end{minipage}} \end{center}

\begin{remark}
The algorithm (\ref{DP-A}) essentially shares the same features as the algorithm (\ref{PD-2}), despite their only difference in the order of updating the primal and dual variables, as well as their very slight difference of the matrices in their respective correction steps.
\end{remark}

\subsection{Specification of the prototypical algorithmic framework (\ref{M-PRECOR})}

Similarly, we analyze the prediction step (\ref{DP-Pre}) and the correction step (\ref{DP-COR}), and show that they can be obtained by specifying the prototypical prediction and correction steps (\ref{M-PRE}) and (\ref{M-COR}), respectively. Hence, the algorithm (\ref{DP-A}) can also be specified by the prototypical algorithmic framework (\ref{M-PRECOR}). The specified matrices are denoted by ${\Q}_{\dup}$ and $ {\M}_{\dup}$, respectively.

\subsubsection{Analysis for the prediction (\ref{DP-Pre})}   \label{Sec41-DP}

According to Lemma \ref{CP-TF},  the predictor generated by \eqref{DP-Pre}  satisfies $\tilde{w}^k\in \Omega$ and
$$
\left\{ \begin{array}{l}         \theta_1(x) - \theta_1(\tilde{x}^k)  + (x - \tilde{x}^k)^T\{-A^T\tilde{\lambda}^k  +\beta A^T\!A(\tilde{x}^k-x^k)\} \ge 0,  \;\; \forall \, x\in {\cal X},\\[0.1cm]
       \theta_2(y) - \theta_2(\tilde{y}^k)  + (y - \tilde{y}^k)^T\{-B^T\tilde{\lambda}^k
           +\beta B^T\!A(\tilde{x}^k-x^k)+ \beta B^T\!B(\tilde{y}^k-y^k) \} \ge 0,  \;\;  \forall \,  y\in {\cal Y}, \\[0.1cm]
  \hskip 1.6cm
   (\lambda - \tilde{\lambda}^k)^T\{\frac{1}{\beta} (\tilde{\lambda}^k-\lambda^k) + (A{x}^k + B{y}^k-b)
              \}\ge 0, \;\;  \forall \, \lambda\in {\Lambda}.
              \end{array}\right.
   $$
Using the VI form \eqref{VI-ID}, we have that
  \[ \label{DP-C}
      \left\{ \begin{array}{l}
       \theta_1(x) - \theta_1(\tilde{x}^k)  + (x - \tilde{x}^k)^T\{\underline{-A^T\tilde{\lambda}^k}  + \beta A^T\!A(\tilde{x}^k-x^k) \} \ge 0,  \quad \forall \, x\in {\cal X},\\[0.1cm]
     \theta_2(y) - \theta_2(\tilde{y}^k)  + (y - \tilde{y}^k)^T\{\underline{-B^T\tilde{\lambda}^k}
                                    +  \beta B^T\!A(\tilde{x}^k-x^k)      \\
 \hskip 6.1cm +\beta B^T\!B(\tilde{y}^k-y^k) \} \ge 0, \quad \forall \, y\in {\cal Y}, \\[0.1cm]
  \hskip 1.0cm
   (\lambda - \tilde{\lambda}^k)^T\{(\underline{A\tilde{x}^k + B\tilde{y}^k-b})  -A(\tilde{x}^k - x^k)
                     -B(\tilde{y}^k - y^k)  \\
            \hskip 6.2cm  + (1/\beta)  (\tilde{\lambda}^k-\lambda^k)\}\ge 0, \quad  \forall \, \lambda\in {\Lambda}.
  \end{array} \right.
              \]
The sum of the underling parts of \eqref{DP-C} is exactly $F(\tilde{w}^k)$, where $F(\cdot)$ is defined in (\ref{VI-Iw}). Thus, we have the following lemma.

\begin{lemma}
\begin{subequations}  \label{DP-PreD}  With the  given
 $(Ax^k, By^k, \lambda^k)$,  the predictor $\tilde{w}^k\in \Omega$ produced by  \eqref{DP-Pre} satisfies
\[ \label{DP-PreD-F }
  \tilde{w}^k\in \Omega, \;\;\;   \theta(u) - \theta(\tilde{u}^k) +  (w- \tilde{w}^k )^T
     \{  F(\tilde{w}^k)  +  Q_{\dup}(\tilde{w}^k-w^k)\}\ge 0,  \quad \forall \, w \in\Omega,\]
where
     \[ \label{DP-Q}
      Q_{\dup}=  \left(\begin{array}{ccc}
          \beta A^TA  &   0                  &        0       \\
          \beta B^TA  & \beta B^TB &         0         \\
                        -A     &   -B    & \frac{1}{\beta} I_m
   \end{array} \right).
    \]
    \end{subequations}
   \end{lemma}
   \noindent{\bf Proof}. The assertion directly follows from  \eqref{DP-C}. \qquad {$\Box$}

   \smallskip

Using the notation in \eqref{XiPw}, we can rewrite the matrix  $Q_{\dup}$ in  \eqref{DP-Q} as
\[  \label{CalQ-DP}   Q_{\dup}= P^T {\Q}_{\dup} P
  \qquad \hbox{where} \qquad  {\Q}_{\dup} = \left(\begin{array}{ccc}
             I_m &   0       &  0  \\
              I_m   &   I_m    & 0   \\
            -I_m  &    -I_m   & I_m
   \end{array} \right).
  \]
Also, note that the matrices  $Q_{\dup}$ and ${\Q}_{\dup}$ are different. It follows from \eqref{DP-PreD} that
\[ \label{DP-PreF}
  \tilde{w}^k\in \Omega, \;\;\;   \theta(u) - \theta(\tilde{u}^k) +  (w- \tilde{w}^k )^T F(\tilde{w}^k)
        \ge (\xi -\tilde{\xi}^k)^T{\Q}_{\dup}(\xi^k - \tilde{\xi}^k),  \quad \forall \, w\in\Omega.\]
Thus, the prediction step (\ref{DP-Pre}) can be specified by the prototypical prediction step (\ref{M-PRE}) with $\Q:= {\Q}_{\dup}$ as defined in (\ref{CalQ-DP}).

\subsubsection{Analysis for the correction procedure (\ref{DP-COR})} \label{Sec42-DP}

Left-multiplying the matrix $\hbox{diag}(\sqrt{\beta}I_m,\sqrt{\beta}I_m, \frac{1}{\sqrt{\beta}}I_m)$  to both sides of the correction step \eqref{DP-COR},  we get
$$   % \label{DP-COR-a}
\left(\begin{array}{c}
   \sqrt{\beta} Ax^{k+1} \\[0.1cm]
  \sqrt{\beta} By^{k+1} \\[0.1cm]
    \frac{1}{\sqrt{\beta}}\lambda^{k+1}
\end{array}\right)   = \left(\begin{array}{c}
\sqrt{\beta}Ax^{k} \\[0.1cm]
\sqrt{\beta}By^{k} \\[0.1cm]
      \frac{1}{\sqrt{\beta}} \lambda^{k}
\end{array}\right)  -\left(\begin{array}{ccc}
         \sqrt{\beta}   {\nu}I_m  &   -  \sqrt{\beta}{\nu} I_m &  0   \\[0.1cm]
                 0 &  \sqrt{\beta} {\nu} I_m & 0  \\[0.1cm]
                - \nu\sqrt{\beta} I_m &    - \nu\sqrt{\beta} I_m  &   \frac{1}{\sqrt{\beta}}  I_m
        \end{array}\right) \left(\begin{array}{c}
Ax^{k} -A\tilde{x}^{k}  \\[0.1cm]
By^{k} - B \tilde{y}^{k}\\[0.1cm]
\lambda^{k} -\tilde{\lambda}^k
\end{array}\right).
  $$
Then, we have
\[ \label{DP-COR-B}
\left(\begin{array}{c}
   \sqrt{\beta} Ax^{k+1} \\[0.1cm]
  \sqrt{\beta} By^{k+1} \\[0.1cm]
    \frac{1}{\sqrt{\beta}}\lambda^{k+1}
\end{array}\right)   = \left(\begin{array}{c}
\sqrt{\beta}Ax^{k} \\[0.1cm]
\sqrt{\beta}By^{k} \\[0.1cm]
      \frac{1}{\sqrt{\beta}} \lambda^{k}
\end{array}\right)  -\left(\begin{array}{ccc}
            {\nu}I_m  &   -  {\nu} I_m &  0   \\[0.1cm]
                 0 &   {\nu} I_m & 0  \\[0.1cm]
                - \nu I_m &    - \nu I_m   &    I_m
        \end{array}\right) \left(\begin{array}{c}
  \sqrt{\beta}(Ax^{k} -A\tilde{x}^{k})  \\[0.1cm]
 \sqrt{\beta} (By^{k} - B \tilde{y}^{k})\\[0.1cm]
    \frac{1}{\sqrt{\beta}} (\lambda^{k} -\tilde{\lambda}^k)
\end{array}\right).
\]
It follows from \eqref{XiPw} that \eqref{DP-COR-B} can be written as
\begin{subequations}  \label{DP-CorB}
\[  \label{DP-CorBw}    {\xi}^{k+1} = {\xi}^k -  {\M}_{\dup}(\xi^k - \tilde{\xi}^k),  \]
where
\[  \label{M-I-DP}  {\M}_{\dup} =
           \left(\begin{array}{ccc}
            {\nu}I_m  &     -{\nu} I_m &  0   \\[0.2cm]
                 0 &   {\nu} I_m & 0  \\[0.2cm]
                -I_m &   -I_m   &   I_m
        \end{array}\right) \quad \hbox{with} \quad \nu\in (0,1).
    \]
\end{subequations}
Thus, the correction step (\ref{DP-COR}) can be specified by the prototypical correction step (\ref{M-COR}) with ${\M}:= {\M}_{\dup}$ as defined in (\ref{M-I-DP}).

\subsection{Convergence}

Since it is shown that the algorithm (\ref{DP-A}) can be specified by the prototypical algorithmic framework (\ref{M-PRECOR}), its convergence can be guaranteed if the conditions (\ref{M-HMQ}) and (\ref{M-HMG}) are satisfied by the just specified matrices ${\Q}_{\dup}$ and ${\M}_{\dup}$. Let us identify the corresponding matrix ${\H}_{\dup}$ with the specified matrices ${\Q}_{\dup}$ and ${\M}_{\dup}$.

\begin{lemma} \label{DP-H} For any $\nu\in(0,1)$, the matrix
\[    \label{MatrixH-DP} {\H}_{\dup} = \left(\begin{array}{ccc}
          \frac{1}{\nu} I_m  &    \frac{1}{\nu} I_m   & 0  \\[0.2cm]
            \frac{1}{\nu} I_m   &   \frac{2}{\nu}   I_m   &  0  \\[0.2cm]
              0   &   0   &   I_m
        \end{array}\right)  \]
 is positive definite, and it holds that
\[   \label{DP-HMQ}   {\H}_{\dup} {\M}_{\dup}= {\Q}_{\dup} , \]
where  ${\M}_{\dup}$ and $ {\Q}_{\dup}$  are defined in \eqref{M-I-DP} and \eqref{CalQ-DP}, respectively.
\end{lemma}
\noindent{\bf Proof}.  It is easy to see that ${\H}_{\dup}$ is positive definite. In addition, we have
\begin{eqnarray*}
{\H}_{\dup} {\M}_{\dup}
     &=&  \left(\begin{array}{ccc}
          \frac{1}{\nu} I_m  &    \frac{1}{\nu} I_m   & 0  \\[0.2cm]
            \frac{1}{\nu} I_m   &   \frac{2}{\nu}   I_m   &  0  \\[0.2cm]
              0   &   0   &   I_m
        \end{array}\right) \left(\begin{array}{ccc}
            {\nu}I_m  &     -{\nu} I_m &  0   \\[0.2cm]
                 0 &   {\nu} I_m & 0  \\[0.2cm]
                -I_m &   -I_m   &   I_m
        \end{array}\right)  \nn   \\
        & =&    \left(\begin{array}{ccc}
             I_m &   0       &   0  \\
              I_m   &   I_m    & 0   \\
             -I_m  &    -I_m   & I_m
   \end{array} \right)     \;=\;  {\Q}_{\dup}.
    \end{eqnarray*}
The assertion \eqref{DP-HMQ} is proved. \qquad {$\Box$}

\begin{lemma}  \label{DP-MG} Let $ {\Q}_{\dup}$, ${\M}_{\dup}$ and ${\H}_{\dup} $  be the matrices defined in \eqref{CalQ-DP}, \eqref{M-I-DP} and \eqref{MatrixH-DP}, respectively. Then the matrix
 \[ \label{DP-G}     {\G}_{\dup}:=(  {\Q}_{\dup}^T +  {\Q}_{\dup})   -
      {\M}_{\dup}^T{\H}_{\dup}{\M}_{\dup}
  \]
is positive definite.
\end{lemma}
\noindent{\bf Proof}. First, by elementary matrix multiplications, we know that
\begin{eqnarray*}
   {\M}_{\dup}^T{\H}_{\dup}{\M}_{\dup}     &  = &  {\M}_{\dup}^T  {\Q}_{\dup}  =
        {\Q}_{\dup}^T  {\M}_{\dup}  \nn \\
        & =&  \left(\begin{array}{ccc}
             I_m &    I_m     &   -I_m \\
              0    &    I_m      &   -I_m  \\
              0  &     0  & I_m
   \end{array} \right)   \left(\begin{array}{ccc}
            {\nu}I_m  &     -{\nu} I_m &  0   \\[0.1cm]
                 0 &   {\nu} I_m & 0  \\[0.1cm]
                -I_m &   -I_m   &   I_m
        \end{array}\right)  \nn \\
        & = & \left(\begin{array}{ccc}
             (1+{\nu})I_m  &    I_m  &  -I_m  \\[0.1cm]
             I_m  &(1 + {\nu}) I_m &   -I_m \\[0.1cm]
                -I_m &   -I_m  &   I_m
        \end{array}\right).
\end{eqnarray*}
Then, we have
\begin{eqnarray*}
   {\G}_{\dup}  &= & (  {\Q}_{\dup}^T +  {\Q}_{\dup})   -
      {\M}_{\dup}^T{\H}_{\dup}{\M}_{\dup}     \nn \\
        & = &  \left(\begin{array}{ccc}
            2 I_m &    I_m     & -  I_m \\
              I_m    &   2 I_m      &  - I_m \\
        -    I_m   &    -I_m   &  2I_m
   \end{array} \right)     -    \left(\begin{array}{ccc}
             (1+{\nu})I_m  &    I_m  &  -I_m  \\[0.1cm]
             I_m  &(1 + {\nu}) I_m &   -I_m \\[0.1cm]
                -I_m &   -I_m  &   I_m
        \end{array}\right)   \nn  \\
   & =&  \left(\begin{array}{ccc}
            (1-{\nu})I_m  &     0 &  0   \\[0.2cm]
                 0 &   (1-{\nu}) I_m &  0  \\[0.2cm]
                  0 &   0  &   I_m
        \end{array}\right).
\end{eqnarray*}
Thus, the matrix  $ {\G}_{\dup}$ is positive definite for any $\nu\in(0,1)$.    \qquad {$\Box$}

\medskip

For the matrices ${\Q}_{\dup}$ in \eqref{CalQ-DP}  and ${\M}_{\dup}$  in \eqref{DP-CorB}, the assertions of Lemmas \ref{DP-H} and \ref{DP-MG} hold. Consequently, we have the following theorem which essentially implies the convergence of the algorithm (\ref{DP-A}). Its proof follows directly from Theorems \ref{HauptA} and \ref{HauptB}.

\begin{theorem} Let $\{ \xi^k\}$ be the sequence generated by the algorithm (\ref{DP-A}). Then, we have
   \[  \label{Haupt-DP}  \|{\xi}^{k+1} -{\xi}^*\|_{ {\H}_{\dup}}^2
       \le   \|{\xi}^k -{\xi}^*\|_{ {\H}_{\dup}}^2  -   \|{\xi}^k - \tilde{\xi}^k  \|_{{\G}_{\dup}}^2,
           \quad \forall \, \xi^*\in {\Xi^*},
       \]
      and thus the sequence $\{ \xi^k\}$  converges to  some  $\xi^* \in \Xi^*$.
 \end{theorem}

\section{Extensions to multi-block separable convex optimization problems with linear equality or inequality constraints}  \label{Sec5-MultiB}

\setcounter{equation}{0}

Recently, there are many intensive discussions on how to extend the ADMM (\ref{Eq-ADMM}) from the two-block separable convex optimization model (\ref{Problem-EqC}) to its generalized  multiple-block models, from both theoretical and algorithmic perspectives. We refer to, e.g., \cite{HTY-SIOPT,HTY-MOR,HY-COAP}, for a few works. Because of the work \cite{CHYY}, it is known that the convergence is not guaranteed if the ADMM (\ref{Eq-ADMM}) is directly extended. This means we should be cautious in both algorithmic design and convergence analysis when considering extensions of the ADMM (\ref{Eq-ADMM}) to multiple-block separable convex optimization problems. As mentioned, this fact also discourages us to reformulate a separable model with linear inequality constraints as another separable model with only linear equality constraints but with more blocks of auxiliary variables, and then consider applying some existing ADMM type algorithms that are eligible to models with linear equality constraints.

In the following three sections, we consider natural extensions from the model (\ref{Problem-IneC}) to its multiple-block generalized one:
\[  \label{Problem-m}
  \min \Bigl\{ \sum_{i=1}^{p} \theta_i(x_i)   \;\big|\;   \sum_{i=1}^{p} A_ix_i=b\ (\hbox{or} \ge b) ,  \;\;  x_i\in {\cal X}_i \Bigr\},\]
where $\theta_i: {\Re}^{n_i}\to {\Re}, \, i=1,\ldots, p$, are closed
proper convex functions and they are not necessarily smooth; ${\cal
X}_i\subseteq \Re^{n_i},\, i=1,\ldots, p$, are closed convex sets;
$A_i\in \Re^{m\times n_i},\, i=1,\ldots, p$, are given matrices; $b\in \Re^m$ is a given vector. The model (\ref{Problem-IneC}) can be regarded as a special case of (\ref{Problem-m}) with $p=2$. Let us focus on the multiple-block case of (\ref{Problem-m}) with $p\ge 3$, and parallelize the discussions in Sections \ref{sec:algframe}-\ref{Sec4-DP} for this multiple-block case. Unlike the failure of the direct extension of the ADMM (\ref{Eq-ADMM}), we will show that the proposed algorithms (\ref{PD-2}) and (\ref{DP-A}) can both be directly extended from the two-block model (\ref{Problem-IneC}) to the multiple-block generalized model (\ref{Problem-m}).

We first summarize some notations and results similar as those in Sections \ref{sec:algframe}-\ref{Sec4-DP}  for the multiple-block model (\ref{Problem-m}). Without ambiguity, some notations are denoted by the same letters as previous sections.

\subsection{VI characterization}

Let $\lambda \in \Re^m$ be the Lagrange multiplier of (\ref{Problem-m}) and the Lagrangian  function  of the problem \eqref{Problem-m} be
\[ \label{Lagrange-F}
  L(x_1,\ldots,x_p,\lambda) =  \sum_{i=1}^{p}\theta_i(x_i) -\lambda^T\Bigl(\sum_{i=1}^{p} A_ix_i-b\Bigr).\]
The optimality condition of (\ref{Problem-m}) can be written as the following VI:
\begin{subequations} \label{VI-FORM}
\[  \label{VI-FORM-Q}
    w^*\in\Omega, \quad \theta(x) - \theta(x^*) + (w-w^*)^T F(w^*) \ge0, \quad
\forall\, w\in \Omega,
    \]
where
\[  \label{VI-FORM-F}    w=\left(\!\!\begin{array}{c}
             x_1\\
         \vdots \\
            x_{p} \\[0.1cm]
           \lambda
             \end{array}\!\! \right), \quad
   x=\left(\!\!\begin{array}{c}
             x_1\\
         \vdots \\
            x_{p} \\
             \end{array}\!\! \right),
    \quad   \theta(x) = \sum_{i=1}^{p} \theta_i(x_i), \quad
       F(w) = \left(\!\!\begin{array}{c}
    - A_1^T\lambda \\
       \vdots \\
    -A_{p}^T\lambda \\[0.1cm]
       \sum_{i=1}^{p} A_ix_i-b
    \end{array}\!\! \right),
  \]
  and
  $$      \Omega = \prod_{i=1}^p {\cal X}_i \times \Lambda \quad \hbox{with}\quad   \Lambda  =\left\{ \begin{array}{ll}
               \Re^m,           &   \hbox{if   $\sum_{i=1}^{p} A_ix_i =  b$} ,   \\[0.2cm]
                \Re^m_+,    &        \hbox{if   $\sum_{i=1}^{p} A_ix_i\ge b$}.
                \end{array}        \right.  $$
\end{subequations}
Again, we denote by $\Omega^*$ the solution set of the VI \eqref{VI-FORM}.

\subsection{Prototypical algorithm framework for VI (\ref{VI-FORM})}\label{sec 6:algframe-1}

Similar as the VI $(\ref{VI-ID})$, we also present a prototypical algorithmic framework for the VI (\ref{VI-FORM}), from which concrete algorithms for the model (\ref{Problem-m}) can be specified. Let us further denote the following notations:

\[ \label{xPD-Bw}
   P =
       \left(\begin{array}{ccccc}
         \sqrt{\beta} A_1    &   \qquad 0 \qquad & \qquad \cdots \qquad & \qquad \cdots \qquad  &    0   \\[0.1cm]
         0   &   \sqrt{\beta} A_2   &  \ddots   &     &    \vdots \\[0.1cm]
             \vdots    &  \ddots   &   \ddots          &    \ddots             &    \vdots \\[0.1cm]
          \vdots   &            &  \ddots  &\sqrt{\beta} A_{p}   &   0\\[0.1cm]
              0   &    \cdots &  \cdots &    0   &  \frac{1}{\sqrt{\beta}}I_{m}
 \end{array}\!\!\right),   \qquad   \xi= Pw= \left(\begin{array}{c}
          \sqrt{\beta} A_1x_1  \\[0.1cm]
             \sqrt{\beta} A_2x_2   \\[0.1cm]
             \vdots \\[0.1cm]
               \sqrt{\beta} A_px_p \\[0.1cm]
            \frac{1}{\sqrt{\beta}} \lambda
   \end{array} \right).\]
Accordingly, we define
$${\Xi}= \big\{\xi  \;|\;  \xi= Pw, \; w\in \Omega\big\}\quad   \hbox{and}  \quad   {\Xi^*}= \big\{\xi^*  \;|\;  \xi^*= Pw^*, \; w^*\in \Omega^*\big\}. $$
Then, the prototypical algorithm framework for the VI (\ref{VI-FORM}) is presented as follows.

 \begin{center}\fbox{
 \begin{minipage}{15.5cm}  {\bf{A Prototypical Algorithmic Framework for VI (\ref{VI-FORM})}}.

\begin{subequations}\label{m-PRECOR}
  \begin{enumerate}
   \item  (Prediction Step) With given $w^k$ and thus $\xi^k= Pw^k$, find $\tilde{w}^k \in \Omega$ such that
  \[  \label{m-PRE} \tilde{w}^k \in \Omega, \;\;\theta(x) - \theta(\tilde{x}^k)  +  (w - \tilde{w}^k)^T F(\tilde{w}^k) \ge
      (\xi  -\tilde{\xi}^k)^T{\Q}(\xi^k-\tilde{\xi}^k), \;\; \forall  \, w \in {\Omega},
        \]
  with $\Q \in \Re^{(p+1)m\times (p+1)m}$, and the matrix ${\Q}^T+{\Q}$  is positive definite.

\item (Correction Step) With $\tilde{w}^k$ solved by (\ref{m-PRE}) and thus $\tilde{\xi}^k=P\tilde{w}^k$, generate $\xi^{k+1}$ by
\[ \label{m-COR}   {\xi}^{k+1} = {\xi}^k -  \M(\xi^k - \tilde{\xi}^k), \]
where $\M \in \Re^{(p+1)m\times (p+1)m}$ is a non-singular matrix.
\end{enumerate}
\end{subequations}
\end{minipage}} \end{center}

\subsection{Roadmap for convergence analysis} \label{sec 6:algframe-2}

Similar as Section \ref{sec:algframe-2}, we prove the convergence of prototype algorithmic framework \eqref{m-PRECOR} for the VI (\ref{VI-FORM}) and set up a roadmap for the convergence analysis.

\begin{theorem} \label{mHauptA}
 For the matrices  ${\Q}$ in \eqref{m-PRE} and $\M$ in \eqref{m-COR},  if  there is a positive definite matrix  $\H\in \Re^{(p+1)m\times (p+1)m}$  such that
  \[ \label{m-HMQ}     \H\M={\Q}  \]
  and
       \[ \label{m-HMG}   \G: = {\Q}^T  + {\Q} - \M^T\H\M\succ 0,  \]
   then we have
  \[  \label{mHauptA0}   \|{\xi}^{k+1} -{\xi}^*\|_{\H}^2
       \le   \|{\xi}^k -{\xi}^*\|_{\H}^2  -   \|{\xi}^k - \tilde{\xi}^k  \|_{\G}^2,
           \quad \forall \, \xi^*\in {\Xi^*}.  \]
   \end{theorem}

Then, analogous to the analysis in Section \ref{sec:algframe-2}, convergence of the prototype algorithmic framework \eqref{m-PRECOR} can be established easily. We summarize the convergence result in the following theorem, and skip the proof.

\begin{theorem} \label{mHauptB}  Let $\{\xi^k\}$ be the sequence generated by the prototype algorithmic framework \eqref{m-PRECOR}.
  If  the conditions \eqref{m-HMQ}  and \eqref{m-HMG} are satisfied, then the sequence $\{ \xi^k\}$ converges to  some
  $\xi^{\infty} \in \Xi^*$.
    \end{theorem}

\subsection{Some useful matrices}

   In order to simplify the notations to be used, we define the following $p\times p$ block matrices:
   \[  \label{x-LowL}
         {\cal L} =  \left(\begin{array}{cccc}
             I_m   &      0      & \cdots      &  0           \\[0.1cm]
             I_m   &     I_m  &  \ddots    &  \vdots   \\[0.1cm]
             \vdots   &  &   \ddots          &       0       \\[0.1cm]
             I_m   &     I_m      &  \cdots  &   I_m
         \end{array}\!\!\right)\qquad \;\hbox{and}\;   \qquad    {\cal I} =
  \left(\begin{array}{cccc}
             I_m   &      0      & \cdots      &  0           \\[0.1cm]
                  0   &     I_m  &  \ddots    &  \vdots   \\[0.1cm]
             \vdots   & \ddots &   \ddots          &       0       \\[0.1cm]
               0  &     \cdots     &  0 &   I_m
         \end{array}\!\!\right).
             \]
We also define the $1\times p$  block matrix
\[  \label{x-Row-E}     {\cal E} =  \left(\!\begin{array}{cccc}
                        I_m   &     I_m      &  \cdots  &   I_m
         \end{array}\!\!\right).   \]
Then, it is easy to see the following properties:

\[  \label{LTL}  {\cal L}^{-1} =    \left(\begin{array}{cccc}
             I_m   &      0      & \cdots      &  0           \\[0.1cm]
             -I_m   &     I_m  &  \ddots    &  \vdots   \\[0.1cm]
              0   &  \ddots &   \ddots          &       0       \\[0.1cm]
               0   &      0     &   -I_m  &   I_m
         \end{array}\!\!\right)
 \quad \hbox{and} \quad  {\cal L}^T + {\cal L} = {\cal I}  + {\cal E}^T {\cal E}.      \]
These matrices will be used in our analysis.

%%%%%%%%%%%%%%%%%%%%%%%%%%%%%%%%%%%%%%%%%%%%%%%%
%%%%%%%%%%%%%%%%%%%%%%%%%%%%%%%%%%%%%%%%%%%%%%%%%

\section{Primal-dual extension of the ADMM (\ref{Eq-ADMM}) for (\ref{Problem-m})}  \label{Sec6}

\setcounter{equation}{0}

This section is parallel to Section \ref{Sec3-PD}. We consider a concrete algorithm for the multiple-block model (\ref{Problem-m}), in which the primal variables $x_i$ ($i=1,\ldots,p$) are updated first before the dual variable $\lambda$. We show that it can be specified from the prototype algorithmic framework \eqref{m-PRECOR}. The penalty parameter is still denoted by $\beta>0$.

\subsection{Algorithm}

\begin{center}\fbox{
 \begin{minipage}{15.5cm}  {\bf{ A Primal-Dual Extension of the ADMM (\ref{Eq-ADMM}) for (\ref{Problem-m})}}.

\begin{subequations}\label{PD-m}
   \begin{enumerate}
   \item (Prediction Step) With given $(A_1x_1^{k},A_2x_2^{k},\cdots, A_px_p^{k}, \lambda^{k})$, find $\tilde{w}^k\in \Omega$ via
\begin{equation}\label{xPD-A}
\left\{
\begin{array}{l}
 \tilde{x}_1^k  \in \arg\min \bigl\{ \theta_1(x_1)  -x_1^TA_1^T\lambda^k   +\frac{\beta}{2} \|A_1(x_1-x_1^k)\|^2    \;|\;   x_1\in{\cal X}_1  \bigr\};\\[0.2cm]
 \tilde{x}_2^k \in \arg\min \bigl\{\theta_2(x_2)  -x_2^TA_2^T\lambda^k   +\frac{\beta}{2} \|A_1(\tilde{x}_1^k-x_1^k)
                 + A_2(x_2-x_2^k)\|^2  \;|\;  x_2\in{\cal X}_2  \bigr\};\\[0.1cm]
  \qquad \qquad \vdots \\
  \tilde{x}_i^k\in\arg\min\bigl\{\theta_i(x_i)  -x_i^TA_i^T\lambda^k   +\frac{\beta}{2} \| \sum_{j=1}^{i-1}A_j(\tilde{x}_j^k-x_j^k)
                 + A_i(x_i-x_i^k)\|^2  \;|\;  x_i\in{\cal X}_i  \bigr\};\\[0.1cm]
  \qquad  \qquad \vdots   \\
  \tilde{x}_{p}^k \in \arg\min \bigl\{
    \theta_{p}(x_{p})  -x_{p}^TA_{p}^T\lambda^k   +\frac{\beta}{2} \| \sum_{j=1}^{p-1}A_j(\tilde{x}_j^k-x_j^k)
                 + A_{p}(x_{p}-x_{p}^k)\|^2
  \;|\;  x_p\in{\cal X}_{p}  \bigr\}; \\[0.3cm]
    \tilde{\lambda}^k =\arg\max\bigl\{- \lambda^T\bigl(\sum_{j=1}^{p} A_j\tilde{x}_j^k -b\bigr) -
       \frac{1}{2\beta}\|\lambda-\lambda^k\|^2   \;|\;  \lambda\in \Lambda \bigr\}.
     \end{array}
\right.
\end{equation}

\item (Correction Step)  Correct the predictor $\tilde{w}^k$ solved by (\ref{m-PRE}), and generate the new iterate $(A_1x_1^{k+1}, A_2x_2^{k+1}, \cdots, A_px_p^{k+1}, \lambda^{k+1})$ with $\nu\in (0,1)$ by
\[ \label{xPD-COR}
\left(\!\!\begin{array}{c}
A_1x_1^{k+1} \\[0.1cm]
A_2x_2^{k+1} \\[0.1cm]
     \vdots \\[0.1cm]
  A_px_p^{k+1} \\[0.1cm]
\lambda^{k+1}
\end{array}\!\!\right)   =
\left(\!\!\begin{array}{c}
A_1x_1^{k} \\[0.1cm]
A_2x_2^{k} \\[0.1cm]
     \vdots \\[0.1cm]
  A_px_p^{k} \\[0.1cm]
\lambda^{k}
\end{array}\!\!\right)
-\left(\!\!\begin{array}{ccccc}
             \nu I_m   &   -\nu I_m     &    0    &     \cdots            &    0 \\[0.1cm]
               0  &   \nu  I_m  &  \ddots    &     \ddots  &          \vdots  \\[0.1cm]
             \vdots   &  \ddots  &   \ddots          &       -\nu I_m        &   0\\[0.1cm]
              0   &    \cdots     &  0  &\nu   I_m &      0\\[0.1cm]
                 - \nu\beta I_m  &     0&  \cdots &   0   &   I_m
 \end{array}\!\!\right)
 \left(\!\!\begin{array}{c}
A_1x_1^{k} -A_1\tilde{x}_1^{k}  \\[0.1cm]
A_2x_2^{k} - A_2\tilde{x}_2^{k}\\[0.1cm]
    \vdots \\[0.1cm]
A_px_p^{k} - A_p\tilde{x}_p^{k}\\[0.1cm]
\lambda^{k} -\tilde{\lambda}^k
\end{array}\!\!\right).
\]
\end{enumerate}
\end{subequations}
\end{minipage}} \end{center}

\begin{remark}
The algorithm (\ref{PD-m}) keeps the main features and structures of various ADMM's extensions in the literature for multiple-block separable convex optimization problems with linear equality constraints, see, e.g., \cite{HTY-SIOPT,HTY-MOR,HY-COAP}. The subproblems in the prediction step (\ref{xPD-A}) treat each $\theta_i$ individually; they are of the same form as those in (\ref{PD-A}) or (\ref{Eq-ADMM}). The correction step (\ref{xPD-COR}) also treats $A_ix_i$ and $A_i{\tilde x}_i$, $i=1,\cdots,p$, aggregately. Hence, it also requires ignorable computation with only few floating-point additions.
\end{remark}

\subsection{Specification of the prototypical algorithmic framework (\ref{m-PRECOR})}

Now, we show that the algorithm (\ref{PD-m}) can be obtained by specifying the prototype algorithmic framework (\ref{m-PRECOR}). That is, we identify the specific matrices $\Q$ and $\M$  in (\ref{m-PRE}) and (\ref{m-COR}) such that (\ref{m-PRE}) and (\ref{m-COR}) can be reduced to the prediction step (\ref{xPD-A}) and the correction step (\ref{xPD-COR}), respectively. The specified matrices corresponding to the algorithm (\ref{PD-m}) are denoted by $\Q_{\pd}$ and ${\M}_{\pd}$, respectively. We also divide the discussion into two subsections.

\subsubsection{Analysis for the prediction step (\ref{xPD-A})}

Similar as Section \ref{Sec31-PD}, for the predictor $\tilde{w}^k$ generated by \eqref{xPD-A}, we have
\begin{subequations}  \label{xPD-PreD}
  \[\label{xPD-PreD-F }
    \theta(x) - \theta(\tilde{x}^k) +(w- \tilde{w}^k)
   ^T F(\tilde{w}^k) \ge (w- \tilde{w}^k)^T
    Q_{\pd}(w^k - \tilde{w}^k),
     \quad \forall \, w\in {\Omega},
       \]
where
\[  \label{xPD-Q}
  Q_{\pd} = \left(\begin{array}{ccccc}
         \beta A_1^TA_1   &      0                         & \cdots      &  0           &    A_1^T   \\[0.3cm]
         \beta A_2^TA_1  &   \beta A_2^TA_2 &  \ddots    &  \vdots   &  A_2^T    \\[0.3cm]
             \vdots              &  &   \ddots          &       0        &    \vdots \\[0.3cm]
         \beta A_{p}^TA_1  &   \beta A_{p}^TA_2      &  \cdots  &\beta A_{p}^TA_{p}  &   A_{p}^T\\[0.3cm]
               0  &     0&  \cdots &   0   &  \frac{1}{\beta}I_{m}
 \end{array}\!\!\right). \]
 \end{subequations}

Using the notations $P$ and $\xi$ in \eqref{xPD-Bw},  and the notations ${\cal L}$ and ${\cal E}$ in \eqref{x-LowL} and \eqref{x-Row-E}, we can rewrite the matrix  $Q_{\pd}$ in  \eqref{xPD-Q} as
   \[  \label{xMatrixBQ-PD}   Q_{\pd}= P^T {\Q}_{\pd} P,
  \quad \hbox{where} \quad  {\Q}_{\pd} = \left(\begin{array}{ccccc}
             I_m   &      0      & \cdots      &  0           &     I_m \\[0.1cm]
             I_m   &     I_m  &  \ddots    &  \vdots   &     I_m   \\[0.1cm]
             \vdots   &  &   \ddots          &       0        &    \vdots \\[0.1cm]
             I_m   &     I_m      &  \cdots  &   I_m &      I_m\\[0.1cm]
               0  &     0&  \cdots &   0   &   I_m
 \end{array}\!\!\right) =\left(\begin{array}{cc}
                         {\cal  L} &      {\cal E}^T\\[0.1cm]
               0  &   I_m
 \end{array}\!\!\right).
  \]
Then, it follows from \eqref{xPD-PreD} that we have the following inequality similar as \eqref{m-PRE}:
\[   \label{xPD-wkw}
     \theta(x) - \theta(\tilde{x}^k)  +  (w - \tilde{w}^k)^T F(\tilde{w}^k) \ge
      (\xi  -\tilde{\xi}^k)^T{\Q}_{\pd}(\xi^k-\tilde{\xi}^k), \quad \forall  \, w \in {\Omega}.
  \]
Thus, the prediction step (\ref{xPD-A}) can be specified by the prototypical prediction step (\ref{m-PRE}) with $\Q:= {\Q}_{\pd}$ as defined in (\ref{xMatrixBQ-PD}).

\subsubsection{Analysis for the correction step (\ref{xPD-COR})}
\medskip

Left-multiplying the matrix $\hbox{diag}(\sqrt{\beta}I_m,\ldots,\sqrt{\beta}I_m, \frac{1}{\sqrt{\beta}}I_m)$  to both sides of the correction step \eqref{xPD-COR},  we get
$$
\left(\!\!\!\begin{array}{c}
 \sqrt{\beta}A_1x_1^{k+1} \\[0.1cm]
 \sqrt{\beta}A_2x_2^{k+1} \\[0.1cm]
     \vdots \\[0.1cm]
  \sqrt{\beta}A_px_p^{k+1} \\[0.1cm]
  \frac{1}{\sqrt{\beta}}\lambda^{k+1}
\end{array}\!\!\!\right)   =
\left(\!\!\!\begin{array}{c}
\sqrt{\beta}A_1x_1^{k} \\[0.1cm]
\sqrt{\beta}A_2x_2^{k} \\[0.1cm]
     \vdots \\[0.1cm]
  \sqrt{\beta}A_px_p^{k} \\[0.1cm]
 \frac{1}{\sqrt{\beta}}\lambda^{k}
\end{array}\!\!\!\right)
-\left(\!\!\!\begin{array}{ccccc}
             \nu \sqrt{\beta} I_m   &   -\nu\sqrt{\beta} I_m     &    0    &     \cdots            &    0 \\[0.1cm]
               0  &   \nu\sqrt{\beta}  I_m  &  \ddots    &     \ddots  &          \vdots  \\[0.1cm]
             \vdots   &  \ddots  &   \ddots          &       -\nu\sqrt{\beta} I_m        &   0\\[0.1cm]
              0   &    \cdots     &  0  &\nu \sqrt{\beta}  I_m &      0\\[0.1cm]
                 - \nu\sqrt{\beta} I_m  &     0&  \cdots &   0   & \frac{1}{\sqrt{\beta}}  I_m
 \end{array}\!\!\!\right)
 \left(\!\!\!\begin{array}{c}
A_1x_1^{k} -A_1\tilde{x}_1^{k}  \\[0.2cm]
A_2x_2^{k} - A_2\tilde{x}_2^{k}\\[0.2cm]
    \vdots \\[0.2cm]
A_px_p^{k} - A_p\tilde{x}_p^{k}\\[0.2cm]
\lambda^{k} -\tilde{\lambda}^k
\end{array}\!\!\!\right).
  $$
It can be written as
\[  \label{PD-AxCOR} \left(\!\!\begin{array}{c}
 \sqrt{\beta}A_1x_1^{k+1} \\[0.1cm]
 \sqrt{\beta}A_2x_2^{k+1} \\[0.1cm]
     \vdots \\[0.1cm]
  \sqrt{\beta}A_px_p^{k+1} \\[0.1cm]
  \frac{1}{\sqrt{\beta}}\lambda^{k+1}
\end{array}\!\!\right)   =
\left(\!\!\begin{array}{c}
\sqrt{\beta}A_1x_1^{k} \\[0.1cm]
\sqrt{\beta}A_2x_2^{k} \\[0.1cm]
     \vdots \\[0.1cm]
  \sqrt{\beta}A_px_p^{k} \\[0.1cm]
 \frac{1}{\sqrt{\beta}}\lambda^{k}
\end{array}\!\!\right)
-\left(\begin{array}{ccccc}
             \nu  I_m   &   -\nu I_m     &    0    &     \cdots            &    0 \\[0.1cm]
               0  &   \nu  I_m  &  \ddots    &     \ddots  &          \vdots  \\[0.1cm]
             \vdots   &  \ddots  &   \ddots          &       -\nu I_m        &   0\\[0.1cm]
              0   &    \cdots     &  0  &\nu   I_m &      0\\[0.1cm]
                 - \nu I_m  &     0&  \cdots &   0   &  I_m
 \end{array}\!\!\right)
 \left(\!\!\begin{array}{c}\sqrt{\beta}(
A_1x_1^{k} -A_1\tilde{x}_1^{k})  \\[0.1cm]
 \sqrt{\beta}(A_2x_2^{k} - A_2\tilde{x}_2^{k})\\[0.1cm]
    \vdots \\[0.1cm]
\sqrt{\beta}(A_px_p^{k} - A_p\tilde{x}_p^{k})\\[0.1cm]
\frac{1}{\sqrt{\beta}} (\lambda^{k} -\tilde{\lambda}^k)
\end{array}\!\!\right).
  \]
Recall the respective definitions ${\cal L}$ and ${\cal E}$ in \eqref{x-LowL} and \eqref{x-Row-E}. We have
        $$
           {\cal L}^{-T} =  \left(\begin{array}{cccc}
             I_m   &   - I_m     &    0    &         0 \\[0.1cm]
               0  &     I_m  &  \ddots    &            0  \\[0.1cm]
             \vdots   &  \ddots  &   \ddots          &       - I_m    \\[0.1cm]
              0   &    \cdots     &  0  &  I_m \\[0.1cm]
               \end{array}\!\!\right) \qquad \hbox{and} \qquad   {\cal E} {\cal L}^{-T}
                  = \left(\!\begin{array}{cccc}
                        I_m   &    0      &  \cdots  &  0
         \end{array}\!\!\right).  $$
Thus, using the notations in \eqref{xPD-Bw}, we can rewrite the correction  step \eqref{PD-AxCOR} as
  \begin{subequations} \label{PD-xiCOR}
  \[  \label{xPD-CORw}    {\xi}^{k+1} = {\xi}^k -  {\M}_{\pd}(\xi^k - \tilde{\xi}^k),  \]
where
\[  \label{xM-PD} {\M}_{\pd} =  \left(\begin{array}{cc}
        {\nu}{\cal L}^{-T}  &  0\\[0.2cm]
     -\nu {\cal E} {\cal L}^{-T}  &   I_m
        \end{array}\right).\]
        \end{subequations}
Thus, the correction step (\ref{xPD-COR}) can be specified by the prototypical correction step (\ref{m-COR}) with $\M:= {\M}_{\pd}$ as defined in (\ref{xM-PD}).

\subsection{Convergence}

Then, according to the roadmap presented in Section \ref{sec 6:algframe-2}, proving the convergence of the algorithm (\ref{PD-m}) can be reduced to verifying the conditions (\ref{m-HMQ}) and (\ref{m-HMG}) with the specified matrices $\Q_{\pd}$ and ${\M}_{\pd}$. That is, the remaining task is to find a positive definite matrix  ${\H}_{\pd}$ such that
  $$      {\H}_{\pd}{\M}_{\pd}={\Q}_{\pd}  \qquad \hbox{and}
   \qquad   {\Q}_{\pd}^T  + {\Q}_{\pd} - {\M}_{\pd}^T{\H}_{\pd}{\M}_{\pd} \succ 0, $$
where ${\Q}_{\pd}$ and ${\M}_{\pd}$  are given by \eqref{xMatrixBQ-PD} and \eqref{xM-PD}, respectively.

\begin{lemma}   \label{PD-M-A}
       For the matrices  ${\Q}_{\pd}$ and ${\M}_{\pd}$  given by \eqref{xMatrixBQ-PD} and  \eqref{xM-PD}, respectively, the matrix
  \[    \label{xMatrixH-PD} {\H}_{\pd} = \left(\begin{array}{cc}
        \frac{1}{\nu}{\cal L}{\cal L}^T + {\cal E}^T{\cal E}   &   {\cal E}^T \\[0.2cm]
 {\cal E}  &   I_m
        \end{array}\right) \quad \hbox{with}  \quad \nu\in (0,1) \]
  is positive definite, and it satisfies  ${\H}_{\pd}{\M}_{\pd} = {\Q}_{\pd} $.
 \end{lemma}
 \noindent{\bf Proof}.   It is easy to check the positive definiteness of $ {\H}_{\pd}$. In addition, for the block matrix ${\Q}_{\pd}$ in \eqref{xMatrixBQ-PD}, we have
 \begin{eqnarray*}
  {\H}_{\pd} {\M}_{\pd} &= &  \left(\begin{array}{cc}
        \frac{1}{\nu}{\cal L}{\cal L}^T + {\cal E}^T{\cal E}   &   {\cal E}^T \\[0.2cm]
 {\cal E}  &   I_m
        \end{array}\right) \left(\begin{array}{cc}
        {\nu}{\cal L}^{-T}  &  0\\[0.2cm]
     -\nu {\cal E} {\cal L}^{-T}  &   I_m
        \end{array}\right)  \nn \\
        &= &\left(\begin{array}{cc}
                         {\cal  L} &      {\cal E}^T\\[0.1cm]
               0  &   I_m
 \end{array}\!\!\right)={\Q}_{\pd}.
 \end{eqnarray*}
 The assertions of this lemma are proved. \qquad {$\Box$}

\begin{lemma}   Let  $ {\Q}_{\pd}$,   ${\M}_{\pd}$  and  ${\H}_{\pd}$ be defined in  \eqref{xMatrixBQ-PD},  \eqref{xM-PD} and \eqref{xMatrixH-PD},  respectively. Then the matrix
 \[ \label{xPD-G}  {\G}_{\pd} := ({\Q}_{\pd}^T +  {\Q}_{\pd})   -
      {\M}_{\pd}^T{\H}_{\pd}{\M}_{\pd}   \]
 is positive definite.
 \end{lemma}
 \noindent{\bf Proof}. By elementary matrix multiplications, we know that
 $${\M}_{\pd}^T{\H}_{\pd}{\M}_{\pd}   =  {\Q}_{\pd}^T  {\M}_{\pd}   =
 \left(\begin{array}{cc}
                         {\cal  L}^T &     0\\[0.1cm]
               {\cal E} &   I_m
 \end{array}\!\!\right)    \left(\begin{array}{cc}
        {\nu}{\cal L}^{-T}  &  0\\[0.2cm]
     -\nu {\cal E} {\cal L}^{-T}  &   I_m
        \end{array}\right)  = \left(\begin{array}{cc}
            {\nu}{\cal I}  &    0 \\[0.1cm]
              0  &   I_{m}
        \end{array}\right).$$
Then, it follows from $ {\cal  L}^T +{\cal L} = {\cal I} + {\cal E}^T{\cal E}$ (see \eqref{x-LowL}-\eqref{LTL}) that
\begin{eqnarray*}
   {\G}_{\pd}  &= & (  {\Q}_{\pd}^T +  {\Q}_{\pd})   -
      {\M}_{\pd}^T{\H}_{\pd}{\M}_{\pd}     \nn \\[0.1cm]
        & = &  \left(\begin{array}{cc}
                         {\cal  L}^T +{\cal L}  &      {\cal E}^T\\[0.1cm]
                 {\cal E}  &  2 I_m
 \end{array}\!\!\right)
             -\left(\begin{array}{cc}
            {\nu}{\cal I}  &    0 \\[0.1cm]
              0  &   I_{m}
        \end{array}\right)  =  \left(\begin{array}{cc}
        (1-\nu){\cal I} + {\cal E}^T {\cal E}  &      {\cal E}^T\\[0.1cm]
                 {\cal E}  &  I_m
 \end{array}\!\!\right).
\end{eqnarray*}
Thus, the matrix  $ {\G}_{\pd}$ is positive definite for any $\nu\in(0,1)$.  \qquad {$\Box$}

Then, according to Theorems \ref{mHauptA} and \ref{mHauptB}, the convergence of the algorithm (\ref{PD-m}) can be obtained. We skip the proof for succinctness.

\section{Dual-primal extension of the ADMM (\ref{Eq-ADMM}) for (\ref{Problem-m})}  \label{Sec7}

\setcounter{equation}{0}

Similar as Section \ref{Sec4-DP}, we can also consider an extension of the ADMM (\ref{Eq-ADMM}) which updates the dual variables $\lambda$ first and then updates the primal variables $x_i$, $i=1,\ldots,p$. The resulting algorithm is called a dual-primal extension of the ADMM (\ref{Eq-ADMM}) for (\ref{Problem-m}). We show that it can also be obtained by specifying the prototype algorithmic framework \eqref{m-PRECOR}.

\subsection{Algorithm}

A dual-primal extension of the ADMM (\ref{Eq-ADMM}) for (\ref{Problem-m}) is presented as follows.

\begin{center}\fbox{
 \begin{minipage}{15.5cm}  {\bf{ A Dual-Primal Extension of the ADMM (\ref{Eq-ADMM}) for  (\ref{Problem-m})}}.

\begin{subequations}\label{DP-m}
   \begin{enumerate}
   \item (Prediction Step) With given $(A_1x_1^{k},A_2x_2^{k},\cdots, A_px_p^{k}, \lambda^{k})$, generate $\tilde{w}^k\in \Omega$ via
\begin{equation}\label{yDP-A}
\left\{
\begin{array}{l}
 \tilde{\lambda}^k =\arg\max\bigl\{- \lambda^T\bigl(\sum_{j=1}^{p} A_j{x}_j^k -b\bigr) - \frac{1}{2\beta}\|\lambda-\lambda^k\|^2    \;|\;  \lambda\in \Lambda \bigr\}; \\[0.2cm]
 \tilde{x}_1^k  \in \arg\min \bigl\{ \theta_1(x_1)  -x_1^TA_1^T\tilde{\lambda}^k   +\frac{\beta}{2} \|A_1(x_1-x_1^k)\|^2  \;|\;   x_1\in{\cal X}_1  \bigr\};  \\[0.2cm]
 \tilde{x}_2^k \in \arg\min \bigl\{\theta_2(x_2)  -x_2^TA_2^T\tilde{\lambda}^k   +\frac{\beta}{2} \|A_1(\tilde{x}_1^k-x_1^k) + A_2(x_2-x_2^k)\|^2  \;|\;  x_2\in{\cal X}_2  \bigr\};  \\[0.1cm]
  \qquad \qquad \vdots \\
  \tilde{x}_i^k\in \arg\min\bigl\{\theta_i(x_i)  -x_i^TA_i^T\tilde{\lambda}^k   +\frac{\beta}{2} \| \sum_{j=1}^{i-1}A_j(\tilde{x}_j^k-x_j^k)
                 + A_i(x_i-x_i^k)\|^2  \;|\;  x_i\in{\cal X}_i  \bigr\};\\[0.1cm]
  \qquad  \qquad \vdots   \\
  \tilde{x}_{p}^k \in\arg\min \bigl\{
    \theta_{p}(x_{p})  -x_{p}^TA_{p}^T\tilde{\lambda}^k   +\frac{\beta}{2} \| \sum_{j=1}^{p-1}A_j(\tilde{x}_j^k-x_j^k)
                 + A_{p}(x_{p}-x_{p}^k)\|^2   \;|\;  x_p\in{\cal X}_{p}  \bigr\} .
     \end{array}
\right.
\end{equation}

\item (Correction Step)  Correct the predictor $\tilde{w}^k$ solved by (\ref{m-PRE}), and generate the new iterate $(A_1x_1^{k+1}, A_2x_2^{k+1}, \cdots, A_px_p^{k+1}, \lambda^{k+1})$ with $\nu\in (0,1)$ by
\[ \label{yDP-COR}
\left(\!\begin{array}{c}
A_1x_1^{k+1} \\[0.1cm]
A_2x_2^{k+1} \\[0.1cm]
     \vdots \\[0.1cm]
  A_px_p^{k+1} \\[0.1cm]
\lambda^{k+1}
\end{array}\!\right)   =
\left(\!\begin{array}{c}
A_1x_1^{k} \\[0.1cm]
A_2x_2^{k} \\[0.1cm]
     \vdots \\[0.1cm]
  A_px_p^{k} \\[0.1cm]
\lambda^{k}
\end{array}\!\right)
-\left(\!\!\begin{array}{ccccc}
             \nu I_m   &   -\nu I_m     &    0    &     \cdots            &    0 \\[0.1cm]
               0  &   \nu  I_m  &  \ddots    &     \ddots  &          \vdots  \\[0.1cm]
             \vdots   &  \ddots  &   \ddots          &       -\nu I_m        &   0\\[0.1cm]
              0   &    \cdots     &  0  &\nu   I_m &      0\\[0.1cm]
                 - \beta I_m  &  -\beta I_m &  \cdots &   -\beta I_m  &   I_m
 \end{array}\!\!\right)
 \left(\!\begin{array}{c}
A_1x_1^{k} -A_1\tilde{x}_1^{k}  \\[0.1cm]
A_2x_2^{k} - A_2\tilde{x}_2^{k}\\[0.1cm]
    \vdots \\[0.1cm]
A_px_p^{k} - A_p\tilde{x}_p^{k}\\[0.1cm]
\lambda^{k} -\tilde{\lambda}^k
\end{array}\!\right).
\]
\end{enumerate}
\end{subequations}
\end{minipage}} \end{center}

\begin{remark}
The algorithm (\ref{DP-m}) differs from the algorithm (\ref{PD-m}) slightly in the order of the update of variables. All subproblems in the prediction step (\ref{yDP-A}) are of the same difficulty as those in (\ref{xPD-A}). The correction step (\ref{yDP-COR}) also differs from (\ref{xPD-COR}) slightly in some entries of their corresponding matrices, and it also only requires ignorable computation.
\end{remark}

\subsection{Specification of the prototypical algorithmic framework (\ref{m-PRECOR})}

Now, we show that the algorithm (\ref{DP-m}) can also be obtained by specifying the prototype algorithmic framework (\ref{m-PRECOR}). The specified matrices $\Q$ and $\M$ in (\ref{m-PRE}) and (\ref{m-COR}) are denoted by $\Q_{\dup}$ and ${\M}_{\dup}$, respectively. Again, we divide the discussion into two subsections.

\subsubsection{Analysis for the prediction step (\ref{yDP-A})}

Similar as the analysis in Section \ref{Sec41-DP}, for the predictor $\tilde{w}^k$ generated by \eqref{yDP-A}, we have
\begin{subequations}  \label{yDP-PreD}
  \[\label{yDP-PreD-F }
    \theta(x) - \theta(\tilde{x}^k) +(w- \tilde{w}^k)
   ^T F(\tilde{w}^k) \ge (w- \tilde{w}^k)^T
    Q_{\dup}(w^k - \tilde{w}^k),
     \; \;\forall\; w\in {\Omega},
       \]
where
\[  \label{yDP-Q}
  Q_{\dup} = \left(\begin{array}{ccccc}
         \beta A_1^TA_1   &      0                         & \cdots      &  0           &   0   \\[0.3cm]
         \beta A_2^TA_1  &   \beta A_2^TA_2 &  \ddots    &  \vdots   &  0   \\[0.3cm]
             \vdots              &  &   \ddots          &       0        &    \vdots \\[0.3cm]
         \beta A_{p}^TA_1  &   \beta A_{p}^TA_2      &  \cdots  &\beta A_{p}^TA_{p}  &  0\\[0.3cm]
              - A_1  &   -  A_2&  \cdots &  - A_p   &  \frac{1}{\beta}I_{m}
 \end{array}\!\!\right). \]
 \end{subequations}

Using the  notations $P$ and $\xi$  in \eqref{xPD-Bw},  and the notations ${\cal L}$ and ${\cal E}$ in \eqref{x-LowL} and \eqref{x-Row-E},  we can rewrite the matrix  $Q_{\dup}$ in  \eqref{yDP-Q} as
   \[  \label{yMatrixBQ-DP}   Q_{\dup}= P^T {\Q}_{\dup} P,
  \quad \hbox{where} \quad  {\Q}_{\dup} = \left(\begin{array}{ccccc}
             I_m   &      0      & \cdots      &  0           &     0 \\[0.1cm]
             I_m   &     I_m  &  \ddots    &  \vdots   &     0   \\[0.1cm]
             \vdots   &  &   \ddots          &       0        &    \vdots \\[0.1cm]
             I_m   &     I_m      &  \cdots  &   I_m &      0\\[0.1cm]
             - I_m  &  - I_m &  \cdots &  - I_m  &   I_m
 \end{array}\!\!\right) =\left(\begin{array}{cc}
                         {\cal  L} &      0\\[0.1cm]
                -{\cal E}  &   I_m
 \end{array}\!\!\right).
  \]
Then, it follows from \eqref{yDP-PreD} that we have the following inequality similar as (\ref{m-PRE}):
\[   \label{yDP-wkw}
     \theta(x) - \theta(\tilde{x}^k)  +  (w - \tilde{w}^k)^T F(\tilde{w}^k) \ge
      (\xi  -\tilde{\xi}^k)^T{\Q}_{\dup}(\xi^k-\tilde{\xi}^k), \quad \forall  \, w \in {\Omega}.
  \]
Thus, the prediction step (\ref{yDP-A}) can be specified by the prototypical prediction step (\ref{m-PRE}) with $\Q:= {\Q}_{\dup}$ as defined in (\ref{yMatrixBQ-DP}).

\subsubsection{Analysis for the correction step (\ref{yDP-COR})}

Left-multiplying the matrix $\hbox{diag}(\sqrt{\beta}I_m,\ldots,\sqrt{\beta}I_m, \frac{1}{\sqrt{\beta}}I_m)$  to both sides of the correction step \eqref{yDP-COR},  we get
$$
\left(\!\!\!\begin{array}{c}
 \sqrt{\beta}A_1x_1^{k+1} \\[0.1cm]
 \sqrt{\beta}A_2x_2^{k+1} \\[0.1cm]
     \vdots \\[0.1cm]
  \sqrt{\beta}A_px_p^{k+1} \\[0.1cm]
  \frac{1}{\sqrt{\beta}}\lambda^{k+1}
\end{array}\!\!\!\right)   =
\left(\!\!\!\begin{array}{c}
\sqrt{\beta}A_1x_1^{k} \\[0.1cm]
\sqrt{\beta}A_2x_2^{k} \\[0.1cm]
     \vdots \\[0.1cm]
  \sqrt{\beta}A_px_p^{k} \\[0.1cm]
 \frac{1}{\sqrt{\beta}}\lambda^{k}
\end{array}\!\!\!\right)
-\left(\!\!\!\begin{array}{ccccc}
             \nu \sqrt{\beta} I_m   &   -\nu\sqrt{\beta} I_m     &    0    &     \cdots            &    0 \\[0.1cm]
               0  &   \nu\sqrt{\beta}  I_m  &  \ddots    &     \ddots  &          \vdots  \\[0.1cm]
             \vdots   &  \ddots  &   \ddots          &       -\nu\sqrt{\beta} I_m        &   0\\[0.1cm]
              0   &    \cdots     &  0  &\nu \sqrt{\beta}  I_m &      0\\[0.1cm]
                 - \sqrt{\beta} I_m  &     0&  \cdots &     - \sqrt{\beta} I_m   & \frac{1}{\sqrt{\beta}}  I_m
 \end{array}\!\!\!\right)
 \left(\!\!\!\begin{array}{c}
A_1x_1^{k} -A_1\tilde{x}_1^{k}  \\[0.2cm]
A_2x_2^{k} - A_2\tilde{x}_2^{k}\\[0.2cm]
    \vdots \\[0.2cm]
A_px_p^{k} - A_p\tilde{x}_p^{k}\\[0.2cm]
\lambda^{k} -\tilde{\lambda}^k
\end{array}\!\!\!\right).
  $$
It can be written as
\[  \label{DP-AxCOR} \left(\!\!\begin{array}{c}
 \sqrt{\beta}A_1x_1^{k+1} \\[0.1cm]
 \sqrt{\beta}A_2x_2^{k+1} \\[0.1cm]
     \vdots \\[0.1cm]
  \sqrt{\beta}A_px_p^{k+1} \\[0.1cm]
  \frac{1}{\sqrt{\beta}}\lambda^{k+1}
\end{array}\!\!\right)   =
\left(\!\!\begin{array}{c}
\sqrt{\beta}A_1x_1^{k} \\[0.1cm]
\sqrt{\beta}A_2x_2^{k} \\[0.1cm]
     \vdots \\[0.1cm]
  \sqrt{\beta}A_px_p^{k} \\[0.1cm]
 \frac{1}{\sqrt{\beta}}\lambda^{k}
\end{array}\!\!\right)
-\left(\!\!\begin{array}{ccccc}
             \nu  I_m   &   -\nu I_m     &    0    &     \cdots            &    0 \\[0.1cm]
               0  &   \nu  I_m  &  \ddots    &     \ddots  &          \vdots  \\[0.1cm]
             \vdots   &  \ddots  &   \ddots          &       -\nu I_m        &   0\\[0.1cm]
              0   &    \cdots     &  0  &\nu   I_m &      0\\[0.1cm]
                 -  I_m  &     0&  \cdots &    -I_m   &  I_m
 \end{array}\!\!\right)
 \left(\!\!\begin{array}{c}\sqrt{\beta}(
A_1x_1^{k} -A_1\tilde{x}_1^{k})  \\[0.1cm]
 \sqrt{\beta}(A_2x_2^{k} - A_2\tilde{x}_2^{k})\\[0.1cm]
    \vdots \\[0.1cm]
\sqrt{\beta}(A_px_p^{k} - A_p\tilde{x}_p^{k})\\[0.1cm]
\frac{1}{\sqrt{\beta}} (\lambda^{k} -\tilde{\lambda}^k)
\end{array}\!\!\right).
  \]
It follows from \eqref{x-LowL} and \eqref{x-Row-E} that
   $$ \left(\begin{array}{ccccc}
             \nu I_m   &   -\nu I_m     &    0    &     \cdots            &    0 \\[0.1cm]
               0  &   \nu  I_m  &  \ddots    &     \ddots  &          \vdots  \\[0.1cm]
             \vdots   &  \ddots  &   \ddots          &       -\nu I_m        &   0\\[0.1cm]
              0   &    \cdots     &  0  &\nu   I_m &      0\\[0.1cm]
                 -  I_m  &    \cdots & -I_m &  -I_m  &   I_m
 \end{array}\!\!\right) = \left(\begin{array}{cc}
        {\nu}{\cal L}^{-T}  &  0\\[0.2cm]
     - {\cal E}  &   I_m
        \end{array}\right).$$
Using the notations in \eqref{xPD-Bw},  we can rewrite the correction  \eqref{DP-AxCOR} as
  \begin{subequations} \label{DP-xiCOR}
  \[  \label{yDP-CORw}    {\xi}^{k+1} = {\xi}^k -  {\M}_{\dup}(\xi^k - \tilde{\xi}^k),  \]
where
\[  \label{yM-DP} {\M}_{\dup} =  \left(\begin{array}{cc}
        {\nu}{\cal L}^{-T}  &  0\\[0.2cm]
     - {\cal E} &   I_m
        \end{array}\right).\]
        \end{subequations}
Thus, the correction step (\ref{yDP-COR}) can be specified by the prototypical correction step (\ref{m-COR}) with $\M:= {\M}_{\dup}$ as defined in (\ref{yM-DP}).

\subsection{Convergence}

Also, according to the roadmap presented in Section \ref{sec 6:algframe-2}, proving the convergence of the algorithm (\ref{DP-m}) can be reduced to verifying the conditions (\ref{m-HMQ}) and (\ref{m-HMG}) with the specified matrices ${\Q}_{\dup}$ and ${\M}_{\dup}$.  That is, the remaining task is to find a positive definite matrix  ${\H}_{\dup}$ such that
$$      {\H}_{\dup}{\M}_{\dup}={\Q}_{\dup}  \qquad \hbox{and}
   \qquad   {\Q}_{\dup}^T  + {\Q}_{\dup} - {\M}_{\dup}^T{\H}_{\dup}{\M}_{\dup} \succ 0, $$
where the matrices $ {\Q}_{\dup}$  and  ${\M}_{\dup}$ are given by  \eqref{yMatrixBQ-DP} and \eqref{yM-DP},  respectively.

\begin{lemma}   \label{DP-M-A}
  For  the matrices $ {\Q}_{\dup}$  and  ${\M}_{\dup}$
    given by  \eqref{yMatrixBQ-DP} and \eqref{yM-DP},  respectively, the matrix
  \[    \label{yMatrixH-DP} {\H}_{\dup} = \left(\begin{array}{cc}
        \frac{1}{\nu}{\cal L}{\cal L}^T   & 0  \\[0.2cm]
    0 &   I_m
        \end{array}\right) \quad \hbox{with} \quad \nu\in (0,1)  \]
is positive definite, and it satisfies   ${\H}_{\dup} {\M}_{\dup} = {\Q}_{\dup} $.
\end{lemma}

\noindent{\bf Proof}. It is clear that  ${\H}_{\dup}$ is positive definite. In addition, for the block matrix ${\Q}_{\dup}$ in \eqref{yMatrixBQ-DP}, we have
 \begin{eqnarray*}
  {\H}_{\dup} {\M}_{\dup} &= &  \left(\begin{array}{cc}
        \frac{1}{\nu}{\cal L}{\cal L}^T  &  0 \\[0.2cm]
 0  &   I_m
        \end{array}\right) \left(\begin{array}{cc}
        {\nu}{\cal L}^{-T}  &  0\\[0.2cm]
     - {\cal E}   &   I_m
        \end{array}\right)   \;=\; \left(\begin{array}{cc}
                         {\cal  L} &    0 \\[0.1cm]
                -{\cal E} &   I_m
 \end{array}\!\!\right)={\Q}_{\dup}.
 \end{eqnarray*}
  The assertions of this lemma are proved. \qquad {$\Box$}

\begin{lemma}   Let $ {\Q}_{\dup}$, ${\M}_{\dup}$ and ${\H}_{\dup}$  be
    defined in  \eqref{yMatrixBQ-DP}, \eqref{yM-DP} and \eqref{yMatrixH-DP},  respectively. Then the matrix
 \[ \label{yDP-G}  {\G}_{\dup} := ({\Q}_{\dup}^T +  {\Q}_{\dup})   -
      {\M}_{\dup}^T{\H}_{\dup}{\M}_{\dup}      \]
 is positive definite.
      \end{lemma}
\noindent{\bf Proof}. First, by elementary matrix multiplications, we get
$${\M}_{\dup}^T{\H}_{\dup}{\M}_{\dup}  =  {\Q}_{\dup}^T  {\M}_{\dup}  \nn  =
\left(\begin{array}{cc}
                         {\cal  L}^T &  - {\cal E}^T\\[0.1cm]
                0 &   I_m
 \end{array}\!\!\right)    \left(\begin{array}{cc}
        {\nu}{\cal L}^{-T}  &  0\\[0.2cm]
     -{\cal E}   &   I_m
        \end{array}\right)  = \left(\begin{array}{cc}
            {\nu}{\cal I} + {\cal E}^T{\cal E}  &    -{\cal E}^T \\[0.1cm]
              - {\cal E}  &   I_{m}
        \end{array}\right).$$
Then,  using $ {\cal  L}^T +{\cal L} = {\cal I} + {\cal E}^T{\cal E}$ (see \eqref{x-LowL}-\eqref{LTL}), we have
\begin{eqnarray*}
   {\G}_{\dup}  &= & (  {\Q}_{\dup}^T +  {\Q}_{\dup})   -
      {\M}_{\dup}^T{\H}_{\dup}{\M}_{\dup}     \nn \\[0.1cm]
        & = &  \left(\begin{array}{cc}
                         {\cal  L}^T +{\cal L}  &    -  {\cal E}^T\\[0.1cm]
                - {\cal E}  &  2 I_m
 \end{array}\!\!\right)
             - \left(\begin{array}{cc}
            {\nu}{\cal I} + {\cal E}^T{\cal E}  &    -{\cal E}^T \\[0.1cm]
              - {\cal E}  &   I_{m}
        \end{array}\right)
                =  \left(\begin{array}{cc}
        (1-\nu){\cal I}  &    0\\[0.1cm]
        0  &  I_m
 \end{array}\!\!\right).
\end{eqnarray*}
Thus, the matrix $ {\G}_{\dup}$ is positive definite for any $\nu\in(0,1)$.

Then, according to Theorems \ref{mHauptA} and \ref{mHauptB}, the convergence of the algorithm (\ref{DP-m}) can be obtained. We skip the proof for succinctness.

\section{Panorama}\label{Sec:overview}

The ALM in \cite{Hes,Powell} was proposed for the nonseparable generic convex optimization problem with linear equality constraints which can be regarded as a one-block model, and as mentioned, the original ADMM (\ref{Eq-ADMM}) is an extension of the ALM by splitting the underlying augmented Lagrangian function twice when the two-block separable model (\ref{Problem-EqC}) is considered. However, as proved in \cite{CHYY}, the same extension may fail in guaranteeing the convergence if a multiple-block generalized model of (\ref{Problem-EqC}) is considered (recall (\ref{Eq-ADMM-z})). Hence, it seems impossible to unify the algorithmic design and convergence analysis for the original ALM in \cite{Hes,Powell}, the original ADMM (\ref{Eq-ADMM}) and its direct extensions, when the number of separable blocks of a convex optimization model with linear equality constraints increases from $p=1$, $p=2$, to $p\ge 3$.

On the other hand, it is easy to see that the algorithms (\ref{PD-m}) and (\ref{DP-m}) are direct extensions of the algorithms (\ref{PD-2}) and (\ref{DP-A}), respectively, when the model under discussion is changed from the two-block model (\ref{Problem-IneC}) to its multiple-block generalized model (\ref{Problem-m}). Alternatively, the algorithms (\ref{PD-2}) and (\ref{DP-A}) are just special cases of the algorithms (\ref{PD-m}) and (\ref{DP-m}) with $p=2$, respectively. Hence, the algorithms (\ref{PD-m}) and (\ref{DP-m}) are eligible to the model (\ref{Problem-m}) with different cases of $p\ge 2$. Indeed, we can show that the algorithms (\ref{PD-m}) and (\ref{DP-m}) can also be applied to the following nonseparable generic convex optimization problem with linear equality or inequality constraints:
\[ \label{Problem-IneC-1}
       \min \big\{\theta(x) \;|\;  A x  =b\ (\hbox{or} \ge b) ,  x\in  {\cal X} \big\},
                \]
which can alternatively be regarded as the special case of (\ref{Problem-m}) with $p=1$.

For the model (\ref{Problem-IneC-1}), let us define
   $$    \Lambda  =\left\{ \begin{array}{ll}
               \Re^m,           &   \hbox{if   $Ax = b$} ,   \\[0.1cm]
                \Re^m_+,    &        \hbox{if   $Ax  \ge b$}.
                \end{array}        \right.
     $$
Then, we can propose an algorithmic framework similar as (\ref{M-PRECOR}) and (\ref{m-PRECOR}) for the generic model (\ref{Problem-IneC-1}), as well as a roadmap for convergence analysis similar as those in Sections \ref{sec:algframe-2} and \ref{sec 6:algframe-2}. Specific algorithms can also be obtained similarly as what we have done in Sections \ref{Sec3-PD}, \ref{Sec4-DP},  \ref{Sec6} and \ref{Sec7}. For succinctness, we only present some specific algorithms and skip other details.

%\subsection{Primal-dual extension of the ADMM for (\ref{Problem-IneC-1})} \label{Sec3-PD}

 \begin{center}\fbox{
 \begin{minipage}{15.5cm}  {\bf{ A Primal-Dual Variant of the ALM for (\ref{Problem-IneC-1})}}.

\begin{subequations}\label{PD-1}
   \begin{enumerate}
   \item (Prediction Step) With given $(Ax^k, \lambda^k)$, find $\tilde{w}^k=(\tilde{x}^k,\tilde{\lambda}^k)$ via
  \[   \label{PD-A1}
  \left\{ \begin{array}{l}
    \tilde{x}^k \in \hbox{argmin}\bigl\{   \theta(x) -  x^TA^T{\lambda}^k  +  \frac{1}{2}\beta \|A(x-x^k)\|^2  \;|\; x\in {\cal X}\bigr\},  \\[0.2cm]
        \tilde{\lambda}^k= \arg\!\max  \bigl\{-\lambda^T\bigl(A\tilde{x}^k-b\bigr)  -     \frac{1}{2\beta}\| \lambda -\lambda^k\|^2  \;|\;  \lambda\in {\Lambda}\bigr\}.
  \end{array}  \right.
\]
\item (Correction Step) Correct the predictor $\tilde{w}^k$ solved by (\ref{PD-A1}), and generate the new iterate $(Ax^{k+1},  \lambda^{k+1})$ with $\nu\in (0,1)$ by
\[ \label{PD-COR-1}
\left(\begin{array}{c}
Ax^{k+1} \\[0.1cm]
\lambda^{k+1}
\end{array}\right)   = \left(\begin{array}{c}
Ax^{k} \\[0.1cm]
\lambda^{k}
\end{array}\right)  -\left(\begin{array}{cc}
                       {\nu} I_m & 0  \\[0.1cm]
                - \nu\beta I_m   &   I_m
        \end{array}\right) \left(\begin{array}{c}
Ax^{k} -A\tilde{x}^{k}  \\[0.1cm]
\lambda^{k} -\tilde{\lambda}^k
\end{array}\right).
\]
\end{enumerate}
\end{subequations}
\end{minipage}} \end{center}

%\subsection{Dual-Primal  extension of the ADMM for (\ref{Problem-IneC-1})} \label{Sec3-PD}

 \begin{center}\fbox{
 \begin{minipage}{15.5cm}  {\bf{ A Dual-Primal Variant of the ALM for (\ref{Problem-IneC-1})}}.

\begin{subequations}\label{DP-A1}
   \begin{enumerate}
   \item (Prediction Step) With given $(Ax^k,  \lambda^k)$, find $\tilde{w}^k=(\tilde{x}^k, \tilde{\lambda}^k)$ via
  \[ \label{DP-Pre-1}
\left\{\begin{array}{l}
                         \tilde{\lambda}^k= \arg\!\max \bigl\{  -\lambda^T \bigl(Ax^k +By^k-b\bigr) - \frac{1}{2\beta}\| \lambda -\lambda^k\|^2   \;|\;  \lambda\in {\Lambda}\bigr\} , \\[0.2cm]
    \tilde{x}^k \in \hbox{argmin}\bigl\{   \theta(x) -  x^TA^T\tilde{\lambda}^k
                         +  \frac{1}{2}\beta \|A(x-x^k)\|^2
                          \;|\;  x\in {\cal X}\bigr\}.
  \end{array} \right.
  \]

\item (Correction Step) Correct the predictor $\tilde{w}^k$ generated by (\ref{DP-Pre-1}), and generate the new iterate $(Ax^{k+1}, \lambda^{k+1})$ with $\nu \in (0,1)$ by
\[ \label{DP-COR-1}
\left(\begin{array}{c}
Ax^{k+1} \\[0.1cm]
\lambda^{k+1}
\end{array}\right)   = \left(\begin{array}{c}
Ax^{k} \\[0.1cm]
\lambda^{k}
\end{array}\right)  -\left(\begin{array}{ccc}
                    {\nu} I_m & 0  \\[0.1cm]
                - \beta I_m &     I_m
        \end{array}\right) \left(\begin{array}{c}
Ax^{k} -A\tilde{x}^{k}  \\[0.1cm]
\lambda^{k} -\tilde{\lambda}^k
\end{array}\right).
    \]
\end{enumerate}
\end{subequations}
\end{minipage}} \end{center}

It is easy to see that, when the special case of (\ref{Problem-IneC-1}) with linear equality constraints is considered, the algorithms (\ref{PD-1}) and (\ref{DP-A1}) differ from the classic ALM in \cite{Hes,Powell} only slightly in their prediction steps with some constant vectors and in their correction steps with ignorable computation. They maintain all major features and structures of the ALM, but they can be used for the cases of (\ref{Problem-IneC-1}) with both linear equality and inequality constraints. In addition,  the algorithms (\ref{PD-1}) and (\ref{DP-A1}) can be rendered from the algorithms (\ref{PD-2}) and (\ref{DP-A}), respectively, by removing the $x_2$-subproblems in (\ref{PD-A}) and (\ref{DP-Pre}) as well as the second rows and columns of the matrices in (\ref{PD-COR}) and (\ref{DP-COR}) correspondingly. Thus, they are also included as special cases by the algorithms (\ref{PD-m}) and (\ref{DP-m}) with $p=1$.

In a nutshell, the proposed algorithmic framework and roadmap for convergence analysis are uniformly eligible to the nonseparable generic convex optimization model (\ref{Problem-IneC-1}), the two-block separable model (\ref{Problem-IneC}), and its multiple-block generalized model (\ref{Problem-m}) with an arbitrary $p$. The resulting algorithms maintain the same features and structures from stem to stern for various convex optimization models with different degrees of separability, in which both linear equality and inequality constraints can be included; and the convergence analysis can be unified by a common roadmap. In this sense, our philosophy of algorithmic design and the roadmap for convergence analysis are panoramic and consistent.

\section{Conclusions} \label{Sec-Conclusion}

The classic alternating direction method of multipliers (ADMM) has been widely used for various convex optimization problems with linear equality constraints and two-block separable objective functions without coupled variables. It is known that the ADMM cannot be directly extended to multiple-block (more than two blocks) separable convex optimization problems with linear equality constraints, while it is unknown whether or not it can be extended to two-block or multiple-block separable convex optimization problems with linear inequality constraints. In this paper, we focus on extensions of the ADMM to both two-block and multiple-block separable convex optimization problems with either linear inequality or linear equality constraints, and propose prototypical algorithmic frameworks which can be specified as concrete algorithms for the targeted models. The specified algorithms keep the major structures and features of the original ADMM, and only require very simple additional steps to guarantee the convergence. We also establish standard roadmaps to prove the convergence of the proposed prototypical algorithmic frameworks without any extra conditions. We show that, if we follow the roadmaps to derive the convergence of any algorithm specified from the proposed prototypical algorithmic frameworks, then essentially it only requires to specify two matrices and then to check the positive definiteness of another matrix. Our analysis is comprehensive enough to uniformly cover the nonseparable generic model as well as the two-block and multiple-block separable convex optimization models, in which both the linear equality and linear inequality constraints can be included. Our analysis only uses very elementary mathematics and hence it is understandable for laymen.

Our aim is to study possible extensions of the original ADMM from a high-level and methodological perspective; thus we do not present any experiment results. As mentioned, the proposed prototypical algorithmic frameworks basically maintain all the major structures and features of the original ADMM (\ref{Eq-ADMM}) which account for its versatility and efficiency, while the additional correction steps are extremely simple in computation. It is easy to empirically verify the efficiency of the algorithms specified from the proposed algorithmic frameworks. For instance, we have tested more than ten benchmark application problems in various fields, including the least absolute shrinkage and selection operator \cite{Tlasso}, the $L_1$ regularized logistic regression problem \cite{Ha-regression}, some basic total-variation-based image reconstruction problems in \cite{CP-Acta}, the support vector machine in \cite{F-Supportvector}, the sparse inverse covariance selection model in \cite{B-MS}, as well as a number of basic optimization models in \cite{Boyd} (including linear and quadratic programming problems, and the least absolute deviations problem). These application problems can all be modelled as concrete applications of the model (\ref{Problem-EqC}) and they have been well solved by the original ADMM (\ref{Eq-ADMM}) in the literatures. For comparison purpose, we implemented the original codes provided by the respective authors and kept their respective well-tuned settings, including the values of the penalty parameter $\beta$, for implementing the prediction steps (\ref{PD-A}) and (\ref{DP-Pre}), and then simply set $\nu=0.99$ for the correction steps (\ref{PD-COR}) and (\ref{DP-COR}). It has been affirmatively verified by our experiments that the proposed algorithms (\ref{PD-2}) and (\ref{DP-A}) perform nearly the same as the original ADMM (\ref{Eq-ADMM}). That is, the versatility and efficiency of the original ADMM (\ref{Eq-ADMM}) are completely maintained by the specified algorithms (\ref{PD-2}) and (\ref{DP-A}) if the special model (\ref{Problem-EqC}) is considered. Here, we opt to skip the tedious descriptions of various numerical results for succinctness. The conclusion is that algorithms specified from the proposed prototypical algorithms frameworks are eligible to the more general models (\ref{Problem-IneC}) and (\ref{Problem-m}), while they can work as well as the original ADMM (\ref{Eq-ADMM}) if the special case (\ref{Problem-EqC}) is considered.

We would like to emphasize that we mainly initiate the foundation of algorithmic design and convergence analysis on the ground of the original ADMM, and our target models are the most generic and abstract separable convex optimization models with linear equality or inequality constraints. We do not further discuss how to modify, specify, or generalize an algorithm that can be specified from the proposed prototypical algorithmic frameworks for the sake of better taking advantage of the structures and properties of a specific application. Hence, we do not discuss how to solve the resulting subproblems more efficiently or how to find better step sizes; nor do we investigate sharper convergence results such as worst-case convergence rates in terms of iteration complexity, various asymptotical convergence rates under different conditions, or other more challenging issues under additional assumptions on the objective functions, coefficient matrices, and/or others. When a specific application problem is considered, it is possible to specify the proposed prototypical algorithmic frameworks as more application-tailored algorithms. It is also possible to discuss how to combine other techniques with the prototypical algorithmic frameworks to obtain more attractive numerical schemes; such examples include acceleration schemes, inertial schemes, neural networks,  stochastic/randomized techniques, and so on. All these more detailed discussions are excluded in our discussion for succinctness. Our focus is exclusively the discussion of extensions of the most fundamental ADMM (\ref{Eq-ADMM}) from the canonical two-block model (\ref{Problem-EqC}) to its generalized two-block model (\ref{Problem-IneC}) and multiple-block model (\ref{Problem-m}), which can include both linear equality and inequality constraints.

\bigskip
%%%%%%%%%%%%%%%%%%%%%%%%%%%%%%%%%%%%%%%%%%%%%%%%%%%%%%%%%%%%%%%%%%%%%%%%%%%%%%%%%%%%%%%%%%%%%%%%%%%%%%%%%
{
}


\baselineskip=12pt  \footnotesize

\begin{thebibliography}{10}

\bibitem{B-MS}
Banerjee, O.,  Ghaoui, L.E.,  d'Aspremont, A.:  Model selection through sparse maximum likelihood estimation for multivariate Gaussian or binary data. J. Mach. Learn. Res. \textbf{9}, 485--516 (2008)

\bibitem{Beck}
Beck, A.: First-Order Methods in Optimization,  MOS-SIAM Series on Optimization  (2017)


\bibitem{Boyd}
Boyd, S., Parikh, N., Chu, E., Peleato, B., Eckstein, J.: Distributed optimization and statistical learning via the alternating direction method of multipliers. Found. Trends Mach. Learn. \textbf{3}(1), 1--122 (2010)


%
%\bibitem{CC2010}
%Chambolle, A., Caselles, V., Cremers, D., Novaga, M., Pock, T.:  An introduction to total variation for image analysis. In: Theoretical Foundations and Numerical Methods for Sparse Recovery, vol. 9 of Radon Series on Computational and Applied Mathematics, pp. 263--340. Walter de Gruyter, Berlin (2010)
%
%
%\bibitem{CHPock}
%Chambolle, A., Pock, T.: A first-order primal-dual algorithm for convex problems with applications to imaging. J. Math. Imaging Vis.  \textbf{40}, 120--145 (2011)


\bibitem{CP-Acta}
Chambolle, A., Pock, T.: An introduction to continuous optimization for imaging.  Acta Numer. \textbf{25}, 161--319 (2016)


\bibitem{CHYY}
Chen, C.H., He, B.S., Ye, Y.Y., Yuan, X.M.: The direct extension of ADMM for multi-block convex minimization problems is not necessary convergent. Math. Program. \textbf{155}, 57--79 (2016)


\bibitem{cher2007}
Cherkassky, V., Mulier, F.: Learning from Data: Concepts, Theory, and Methods. Wiley-IEEE Press, New York (2007)


\bibitem{CS2000}
Cristianini, N., Shawe-Taylor, J.: An Introduction to Support Vector Machines. Cambridge University Press, Cambridge (2000)


\bibitem{Eck12}
Eckstein, J., Yao, W.: Augmented Lagrangian and alternating direction methods for convex optimization: A tutorial and some illustrative computational results. Pac. J. Optim. \textbf{11}(4), 619--644 (2015)


\bibitem{F-Supportvector}
 Forero, P.A.,  Cano, A., Giannakis,  G.B.:  Consensus-based distributed support vector machines. J. Mach. Learn. Res. \textbf{11}, 1663--1707 (2010)


\bibitem{Gabay}
 Gabay, D.: Application of the method of multipliers to varuational inequalities. In: Fortin, M., Glowinski, R. (eds.) Augmented Lagrangian Methods: Application to the Numerical Solution of Boundary-Value Problem, pp. 299--331. North-Holland, Amsterdam (1983)


\bibitem{Glow84}
Glowinski, R.: Numerical Methods for Nonlinear  Variational Problems. Springer, Berlin (1984)

\bibitem{Glow12}
Glowinski, R.: On alternating direction methods of multipliers: a historical perspective. In: Fitzgibbon, W., Kuznetsov, Y., Neittaanm$\ddot{\hbox{a}}$ki, P., Pironneau, O. (eds.) Modeling, Simulation and Optimization for Science and Technology. Computational Methods in Applied Sciences, vol. 34. Springer, Dordrecht (2014)


\bibitem{GM}
Glowinski, R., Marrocco, A.: Approximation par $\acute{\hbox{e}}$l$\acute{\hbox{e}}$ments finis d'ordre un et r$\acute{\hbox{e}}$solution par p$\acute{\hbox{e}}$nalisation-dualit$\acute{\hbox{e}}$ d'une classe de probl$\acute{\hbox{e}}$mes non lin$\acute{\hbox{e}}$aires.  RAIRO Anal. Numer. R2, 41--76 (1975)

\bibitem{Glow89}
 Glowinski, R., Le Tallec, P.: Augmented Lagrangian and Operator-Splitting Methods in Nonlinear Mechanics. SIAM Studies in Applied Mathematics, Philadelphia, PA (1989)


\bibitem{Ha-regression}
Hastie, T., Tibshirani,  R., Friedman, J.: The Elements of Statistical Learning: Data Mining, Inference and Prediction.  Springer-Verlag, New York (2017)


\bibitem{He20}
He, B.S.: My 20 years research on alternating directions method of multipliers (Chinese). Oper. Res. Trans.  \textbf{22},  1-31  (2018)


\bibitem{He02MP}
He, B.S., Liao, L.Z.,  Han, D., Yang, H.:  A new inexact alternating directions method for monontone variational inequalities. Math. Program., \textbf{92}, 103--118 (2002)


\bibitem{HMY-COAP}
He, B.S., Ma, F., Yuan, X.M.: Optimally linearizing the alternating direction method of multipliers for convex programming. Comput. Optim. Appl. \textbf{75}(2), 361--388 (2020)


\bibitem{HTY-SIOPT}
He, B.S., Tao, M., Yuan, X.M.: Alternating direction method with Gaussian back substitution for separable convex programming. SIAM J. Optim. \textbf{22}(2), 313--340 (2012)


\bibitem{HTY-MOR}
He, B.S., Tao, M., Yuan, X.M.: Convergence rate analysis for the alternating direction method of multipliers with a substitution procedure for separable convex programming. Math. Oper. Res. \textbf{42}(3), 662--691 (2017)


\bibitem{HY-SINUM}
He, B.S., Yuan, X.M.: On the $\mathcal{O}(1/n)$ convergence rate of Douglas-Rachford alternating direction method. SIAM J. Numer. Anal. \textbf{50}, 700--709 (2012)


\bibitem{HY-NM}
He, B.S., Yuan, X.M.: On non-ergodic convergence rate of Douglas-Rachford alternating directions method of multipliers, Numer. Math., \textbf{130}, 567--577 (2015)


\bibitem{HY-COAP}
He, B.S., Yuan, X.M.: A class of ADMM-based algorithms for three-block separable convex programming. Comput. Optim. Appl. \textbf{70}, 791--826 (2018)


\bibitem{Hes}
Hestenes, M.R.:  Multiplier and gradient methods. J. Optim. Theory Appli. \textbf{4}, 303--320 (1969)


\bibitem{Lee2001}
Lee, Y.J., Mangasarian, O.L.: SSVM: a smooth support vector machines for classification. Comput. Optim. Appl. \textbf{20}(1), 5--22 (2001)



\bibitem{Man2000}
Mangasarian, O.L.: Generalized support vector machines. In: Smola, A., Bartlett, P., Scholkopf, B., Schuurmans, D. (eds) Advances in Large Margin Classifiers,  pp. 135--146, MIT Press, Cambridge  (2000)



\bibitem{Powell}
Powell, M.J.D.: A method for nonlinear constraints in minimization problems. In: Fletcher, R. (ed.) Optimization, pp. 283--298. Academic Press, New York (1969)



%\bibitem{ROF1992}
%Rudin, L.I., Osher, S., Fatemi, E.: Nonlinear total variation based noise removal algorithms. Physica D \textbf{60}(1-4),   259--268 (1992)


\bibitem{Tlasso}
Tibshirani, R.: Regression shrinkage and selection via the lasso. J. R. Stat. Soc. Ser. B \textbf{58}(1), 267--288 (1996)


\bibitem{Vapn1995}
 Vapnik, V.N.: The Nature of Statistical Learning Theory. Springer, New York (1995)


\bibitem{YZZ-JMLR}
Yuan, X., Zeng, S., Zhang, J.: Discerning the linear convergence of ADMM for structured convex optimization through the lens of variational analysis. J. Mach. Learn. Res. \textbf{21}, 1--75 (2020)


\end{thebibliography}
\end{document}